\documentclass[12pt,a4paper]{article}

\usepackage{amsmath,amssymb,amsthm}
\usepackage{array}
\usepackage{mathtools}
\usepackage[utf8]{inputenc}
\usepackage[T1]{fontenc}
\usepackage{physics}
\usepackage[margin=1in]{geometry}
\usepackage{hyperref}
\usepackage{enumitem}
\usepackage{bbm}
\usepackage{braket}
\usepackage{graphicx}
\usepackage[backend=biber, style=numeric]{biblatex}
\usepackage{xcolor}
\usepackage{authblk}   % preamble

\newtheorem{theorem}{Theorem}
\newtheorem{lemma}[theorem]{Lemma}
\newtheorem{corollary}[theorem]{Corollary}
\theoremstyle{definition}
\newtheorem{definition}{Definition}
\theoremstyle{remark}
\newtheorem{remark}{Remark}

\newcommand{\R}{\mathbb{R}}
\newcommand{\hc}{{\mathrm{hc}}}

\newcommand{\Bone}{B_{1}(0)}

\title{The Lee--Huang--Yang Formula\\in the Hard-Core Simplified Approach to the Bose gas}

\author[1]{M.E.H. Bahri\thanks{\texttt{mbahri@math.rutgers.edu}}}
\author[1]{Ian Jauslin\thanks{\texttt{ian.jauslin@rutgers.edu}}}
\affil[1]{Department of Mathematics, Rutgers University, 08854, U.S.A.}

\date{}

\begin{document}
\maketitle

\begin{abstract}
We study the hard-core variant of the Simplified Approach to the ground state of the dilute Bose gas.
The hard-core potential has been the subject of renewed attention, and is physically relevant to realistic atomic Bose gases.
The hard-core Simple Equation is an integro-differential equation, and
%for the function \(\tau=1-g\), where \(g\) is the pair correlation function, so that \(\tau\equiv1\) inside the hard core of radius \(1\).
we construct its solution via a monotone iteration scheme and prove existence and uniqueness (in a restricted sense) of the solution. The energy per particle \(e\) and the density \(\rho\) of the solution are shown to satisfy the Lee--Huang--Yang formula
\[
e = 2\pi\rho\Bigl(1+\frac{128}{15\sqrt{\pi}}\sqrt{\rho}+o(\sqrt{\rho})\Bigr), \qquad \rho\to 0,
\]
and the condensate fraction \(\eta\) predicted by the hard-core Simple Equation is shown to satisfy
\[
\eta = 1-\frac{8}{3\sqrt{\pi}}\sqrt{\rho}+o(\sqrt{\rho}), \qquad \rho\to 0,
\]
in agreement with Bogolyubov's theory.
We thus extend results of Carlen, Jauslin, Lieb to the hard-core case.
The hard-core equation is defined on its own terms, and the proofs proceed by approximating the hard-core potential with the soft-core potentials \(m\mathbbm{1}_{\{|x|\le1\}}\) of increasing height \(m\): at fixed energy, we establish the \(L^1\) convergence of the corresponding solutions and the convergence of their densities at an \(O(m^{-1/2})\) rate, and we show that the Lee--Huang--Yang expansions of the energy and of the condensate fraction for the soft-core equations hold with errors uniform in \(m\), so that they survive the hard-core limit.
%For the condensate fraction this also requires interchanging the limit \(m\to\infty\) with the derivative \(\partial_\mu|_{\mu=0}\) that defines it.
\end{abstract}

%============= INTRODUCTION =============
\section{Introduction}
It has long been known that the Bose gas condenses below a critical temperature \cite{bo1924,e1925} as was verified experimentally by \cite{Ae95,De95} in extremely low density samples of rubidium and sodium. Despite phenomenological theories that were developed by Bogolyubov, Gross, Pitaevskii, Landau and others \cite{g1961,p1961,bog1947,lan1941} in the low-density regime, a proof of Bose-Einstein condensation from first principles in interacting Bose gases is lacking (very recently, in a release of machine-generated manuscripts \cite{OAI}, a proof of Bose--Einstein condensation for the dilute hard-sphere Bose gas in the thermodynamic limit has been announced; this result is still under review, and is yet to be published in a scientific journal).
Proving many other properties of interacting Bose gases including super-fluidity remains an open problem.

Over the past three decades, starting with Lieb and Yngvason's study of its ground state energy \cite{LY98, LSSY05}, mathematical physicists have worked towards filling these gaps. There have been many successes both in the Gross-Pitaevskii regime \cite{LS02,Bo18} and the mean-field regime \cite{S11}, and in the thermodynamic limit \cite{YY09,FS20,FS22,BCG24,BBCOS26,BCS21}. Lately there has been renewed interest in studying the hard-core Bose gas \cite{LY98,FS22,BBCOS26,BCG24}. This is for two reasons. Firstly, it is a mathematical challenge due to the infinite gradient of the potential and secondly realistic potentials for bosonic atoms such as the Aziz potential \cite{ANe79} or the HFDHE2 potential \cite{BC94} include a strongly repulsive singularity at the origin which is well approximated by a hard-core.

Most of the results that have been obtained in the thermodynamic limit concern the ground-state energy \cite{YY09,FS20,FS22,BBCOS26,BCG24} and free energy \cite{BBCD25, FGJMOT26, S08,Y10,He23}. However a description of other observables of Bose gases, such as the condensate fraction, still remains lacking \cite{LSSY05}. 

%The condensate fraction of every ground state was found to be bounded below by an absolute constant once $\rho a^3$ is below an absolute threshold (``Ground-state condensation in the dilute hard-sphere gas'', September 24, 2026). In addition, the condensation persists at a fixed small positive temperature for the canonical Gibbs state at small fixed density (``Bose--Einstein condensation at positive temperature in the dilute hard-sphere gas'', October 5, 2026), and the leading-order Bogolyubov depletion law is obtained for the hard-sphere gas (``Quantum depletion and momentum distribution in the dilute hard-sphere Bose gas'', October 5, 2026). The main ideas of these proofs are credited to ++++++}

In 1963, Lieb introduced an effective equation for the ground state properties of the interacting Bose gas in the thermodynamic limit, known as the Simplified Approach \cite{Li63}. In a recent sequence of papers \cite{CJL20,CJL21,CHJL21,J25}, the Simplified Approach has been shown to be surprisingly accurate. In particular, its prediction for the ground-state energy asymptotically matches that for the Bose gas both at low \emph{and} high densities \cite{CJL20}. In addition, its predictions for the condensate fraction, the momentum distribution and the $2$-point correlation function agree with Bogolyubov's theory at low densities \cite{CJL21}. Furthermore, numerical studies showed the ground-state energy prediction of the Simplified Approach is quantitatively accurate when compared with Monte-Carlo simulations at \emph{all densities} as well as the qualitative behavior of the condensate fraction \cite{CHJL21}.

The existing results of the Simplified Approach have been established for integrable soft-core potentials. In this note, we develop a theory for the Simplified Approach with a hard-core potential. We prove that the ground-state energy with a hard-core potential agrees asymptotically with the many-body Bose gas at low densities (up to Lee-Huang-Yang [LHY] order; see \eqref{eq:lhyeq}). In addition, we show that the condensate fraction is in agreement with Bogolyubov's theory.

The derivation of the Simplified Approach from the many-body Bose gas is non-rigorous and relies on a truncation of a hierarchy of correlation functions in terms of the $2-$point correlation function $g(x)$ \cite{Li63,J25}. This approximation leads to a collection of equations of varying difficulties of complexity and accuracy. On one side of the spectrum the Complete Equation was numerically verified to be accurate as mentioned previously but is exceedingly difficult to study analytically \cite{CHJL21}. On the other end of the spectrum is the Simple Equation which is only accurate at high and low densities, but is simple enough to be studied analytically. In this note, we will prove properties of the Simple Equation with a hard-core potential at low densities (since there is a critical packing fraction the high-density regime is inaccessible). For a numerical analysis of the solutions see \cite{CHJL21}.

\begin{definition}[Simple Equation]
The Simple Equation of the ground-state Bose gas in $\mathbb{R}^3$ with a non-negative radially-symmetric potential, $v(x) \in L^1(\mathbb{R}^3) \cap L^p(\mathbb{R}^3)$ with $p>\frac{3}{2}$ is
\begin{align}
(-\Delta + 4e + v(|x|))u(|x|) & = v(|x|) + 2e \rho u \ast u (|x|), \label{eq:simpleq}\\
e &= \frac{\rho}{2} \int (1-u(|x|))v(|x|)d x \label{eq:simpleq-e}
\end{align}
where $u$ is an integrable and radially-symmetric function. $e$ is the prediction for the ground state energy per particle.
\end{definition}

\cite{CJL20} proved the existence of a solution satisfying $0 \le u \le 1$ and tied the relationship between the density and energy (given via the second equation in \eqref{eq:simpleq}) to the integrability and uniqueness of the solution.

In addition, the energy-per-particle, $e$, of the solution was also found to agree with the seminal Lee-Huang-Yang conjecture at low density in 3 dimensions:
\begin{equation}\label{eq:lhyeq}
e = 2 \pi \rho a \bigg[ 1 + \frac{128}{15 \sqrt{\pi}}\sqrt{\rho a^3} + o(\sqrt{\rho})\bigg]
\end{equation}
where $a$ is the scattering length of the corresponding $2$-body problem with potential $v$. The LHY formula was recently proved for the soft-core many-body Bose gases by \cite{YY09,FS20,BCS21}.

To introduce the hard-core Simple Equation, we note that \eqref{eq:simpleq} is equivalent to
\begin{align}
(-\Delta + 4e + v(|x|))u(|x|) & = v(|x|) + 2e \rho u \ast u (|x|), \\  \int u &= \frac{1}{\rho} \label{eq:simpleq-eq}
\end{align}
Indeed, when the potential is integrable, integration of the first equation of \eqref{eq:simpleq} yields
\begin{equation}
2 e \rho \bigg( \int u \bigg)^2 - 4 e \int u  + \int [1-u]v = 0
\end{equation}
which implies that
\begin{equation}
e = \frac{\rho}{2} \int (1-u)v \Leftrightarrow \int u = \frac{1}{\rho}.
\end{equation}
This re-formulation naturally suggests the generalization of the definition of the Simple Equation to the hard-core potential:
since the wave-function vanishes inside the hard core we impose $1 -\tau(x)=0$ for $|x| \le 1$:
\begin{definition}[Hard-Core Simple Equation]
\begin{align}
\Delta \tau(|x|) &= 4 e \tau(|x|) - 2 e \rho_{\hc} \tau \ast \tau (|x|), \quad &&|x| > 1, \label{eq:hsimpleq}\\
\tau(|x|) &=1, \quad &&|x| \le 1, \\
\int \tau &= \frac{1}{\rho_{\hc}}. \label{eq:hsimpleq-norm}
\end{align}
\end{definition}

We can recover a formula for the ground state energy from (\ref{eq:hsimpleq-norm}) by integrating \eqref{eq:hsimpleq} over $\{|x|>1\}$:
\begin{equation}\label{eq:ehard-pre}
e = \frac{8\pi}3 e \rho_{\hc} - 2\pi\rho_{\hc} \partial_r \tau\big|_{1^+} - e \rho_{\hc}^2 \int_{\Bone} \tau \ast \tau,
\end{equation}
where $\Bone$ is the unit ball centered at the origin, $\frac{4\pi}3 = |\Bone|$ is its volume, and $4\pi$ is the area of the unit sphere.
From (\ref{eq:ehard-pre}), we get
\begin{equation}\label{eq:ehard}
e = \frac{-2\pi\rho_{\hc} \partial_r \tau\big|_{1^+}}{1-\frac{8\pi}3 \rho_{\hc}  + \rho_{\hc}^2 \int_{\Bone} \tau \ast \tau}
.
\end{equation}

Our first main result is the following.

\begin{theorem}[LHY Formula]\label{thm:LHY}
The hard-core Simple Equation \eqref{eq:hsimpleq} has a unique integrable solution $0\le \tau \le 1$ and, as $\rho_{\hc} \rightarrow 0$
\begin{equation}\label{eq:LHY}
	e = 2\pi\rho_{\hc} a_0 \left(1 + \frac{128}{15\sqrt{\pi}}\sqrt{\rho_{\hc} a_0^3} + o\bigl(\sqrt{\rho_{\hc}a_0^3}\bigr)\right),
\end{equation}
where $a_0=1$ is the scattering length of the hard-sphere potential.
\end{theorem}
In other words we recover the LHY formula (\ref{eq:lhyeq}) for the hard-core Simple Equation. The LHY formula was proved for the hard-core many-body Bose gas by \cite{FS22,BBCOS26,BCG24} and thus we show that the energy of the hard-core Simple Equation is asymptotically accurate in the low density limit.

\subsection{Condensate Fraction}\label{sec:observables}
The Simplified Approach also makes predictions about quantities other than the ground state energy \cite{CJL21,J25}; here we focus on the condensate fraction.
Denote by $P_j$ the projector of the $j$-th particle onto the condensate state, which is the constant state:
\begin{equation}
P_j \psi(x_1,\dots,x_N) = \int \frac{d x_j}{V} \psi(x_1,\dots,x_N),
\end{equation}
so that $\frac1N\sum_{j=1}^N P_j$ is the observable for the occupation fraction of the Bose-Einstein condensate.
Its prediction is obtained using the Hellmann-Feynman trick: we modify the many-body Hamiltonian as
\begin{equation}
\mathcal{H}_{N,\mu} := -\frac{1}{2} \sum_{i=1}^N \Delta_i + \sum_{i<j} v(|x_i-x_j|) - \frac{\mu}{N} \sum_{j=1}^N P_j
\label{eq:Hmu}
\end{equation}
and, denoting its ground state by $\psi_{0,\mu}$ and its ground state energy by $\mathcal{E}_{0,\mu}$, we have
\begin{equation}
	 \braket{\psi_0 , \tfrac1N\textstyle\sum_{j=1}^N P_j\, \psi_0}=- \partial_\mu \mathcal{E}_{0,\mu} |_{\mu=0} \label{eq:en}.
\end{equation}
The prediction of the Simplified Approach for the condensate fraction is obtained by deriving the Simplified Approach from the modified Hamiltonian (\ref{eq:Hmu}) (following the same prescription as described above), computing the prediction for its ground state energy $e_\mu$, and differentiating with respect to $\mu$.
Since $P_j$ is a one-particle operator, this can be done in a straightforward way \cite[(1.15)]{CJL21}: the perturbation enters the Simple Equation only through the constant shift $-\Delta\mapsto-\Delta+2\mu$, which yields the modified Simple Equation
\begin{align}
(-\Delta + 2\mu + 4e_\mu + v(|x|))u_\mu(|x|) = & v(|x|) + 2e_\mu \rho u_\mu \ast u_\mu (|x|),\\
e_\mu = & \frac{\rho}{2} \int (1-u_\mu(|x|))v(|x|)d x \label{eq:obssimpleq}
\end{align}
The modified {\it hard-core} Simple Equation is obtained analogously. For $\mu\neq0$, the energy relation \eqref{eq:obssimpleq} is no longer equivalent to $\int u_\mu=1/\rho$: integrating the modified Simple Equation as above shows that it is equivalent to $e_\mu\bigl(1-\rho\int u_\mu\bigr)^2=\mu\,\rho\int u_\mu$, and this is the relation that carries over to the hard core:
\begin{align}
\Delta \tau_\mu(|x|) &= (4 e_\mu + 2\mu) \tau_\mu(|x|) - 2 e_\mu \rho_{ \hc} \tau_\mu \ast \tau_\mu (|x|), \quad &&|x| > 1, \label{eq:obshsimpleq}\\
\tau_\mu(|x|) &=1, \quad &&|x| \le 1, \\
e_\mu\Bigl(1-\rho_{\hc}\int \tau_\mu\Bigr)^2 &= \mu\,\rho_{\hc}\int \tau_\mu, \quad &&\rho_{\hc}\int \tau_\mu\le 1. \label{eq:obshsimpleq-norm}
\end{align}
This yields a prediction for the condensate fraction:
\begin{equation}
  \eta:=-\partial_\mu e_\mu
  =1+  \frac{\rho}{2} \int v(x) \partial_\mu u_\mu (x) |_{\mu=0} d x, \label{eq:condensatefract}
\end{equation}
which \cite{CJL21} proved to scale as
\begin{equation}
\eta = 1- \frac{8 \sqrt{\rho a_0^3}}{3 \sqrt{\pi}} + o(\sqrt{\rho}),
\end{equation}
coinciding with Lee, Huang and Yang's predictions \cite{LHY57}. Our next main result shows that the hard-core equation also possesses the same condensate fraction.

\begin{theorem}[Condensate Fraction]\label{thm:cf}
	The hard-core Simple Equation \eqref{eq:hsimpleq}, in the dilute limit $\rho_{\hc} \rightarrow 0$, satisfies
\begin{equation}
	\eta = 1- \frac{8 \sqrt{\rho_{\hc}\, a_0^3}}{3 \sqrt{\pi}} + o(\sqrt{\rho_{\hc}}), \label{eq:cfeq}
\end{equation}
where $a_0=1$ is the scattering length.
\end{theorem}

\subsection{Proof strategy: limit of strong soft cores}

%{\color{red}Short intro to emphasize that we don't define the hard core as a limit of soft cores, we prove it.}
%To prove the above result and our next result about the expectation of observables, we study the soft-core Simple Equation, variants of it and their iterative solutions as the amplitude of the soft-core potential tends to infinity.
%The hard-core and non-hard-core Simple Equations are independently and heuristically introduced via Lieb's Simplified Approach.Rather than deriving the hard-core as a limiting case of soft-cores, we prove it explicitly in this paper.  

The main tool in the proofs is to compute the energy and the condensate fraction of the hard-core Simple Equation \eqref{eq:hsimpleq} as limits of the energy and the condensate fraction for a {\it soft-core} Simple Equation.
We will use the soft-core potential
\begin{equation}
  v_m(|x|):=m \mathbbm{1}_{\Bone}(x)
  :=
  \left\{\begin{array}{>\displaystyle ll}
    m&\mathrm{if\ }x\in\Bone:=\{x:\ |x|\leqslant 1\}\\
    0&\mathrm{otherwise}
    .
  \end{array}\right.
\end{equation}
The soft-core Simple Equation with potential $v_m$ is then (see \eqref{eq:simpleq}-\eqref{eq:simpleq-e})
\begin{equation}
(-\Delta + 4e + v_m(|x|))u(|x|)  = v_m(|x|) + 2e \bar\rho_m u \ast u (|x|),\quad
e = \frac{\bar\rho_m}{2} \int (1-u(|x|))v_m(|x|)d x
.
\label{eq:ssimpleq}
\end{equation}
We will prove that the energy and the condensate fraction for the hard-core equation are the limit $m\to\infty$ of the soft-core ones. We stress that the hard-core Simple Equation~\eqref{eq:hsimpleq} was introduced above on its own terms, not as a limit of soft-core equations: the two are independently defined objects, each with its own existence and uniqueness theory (see Lemmas \ref{lem:convergence} and \ref{lem:uniqueness} below for the hard-core case).
%The content of Theorem~\ref{thm:OE} is that they are nevertheless compatible: solutions, energies, and observable expectations of the Simple Equation with potential $m\mathbbm{1}_{\Bone}$ converge, at a quantitative rate, to their hard-core counterparts as $m\to\infty$. The limit relation is thus a theorem rather than a convention, which strengthens the robustness of the Simplified Approach: its predictions vary continuously in the interaction, even through the singular hard-core limit.

%\begin{theorem}[Observable Expectation]\label{thm:OE}
%{\color{red}first result: $e=\lim_{m\to\infty}e_m$}
%Consider an observable $\mathcal{U}$ that modifies the Simple Equation per the prescription above with a soft-core potential $m \mathbbm{1}_{\Bone}$, then the following limiting expectations are equal
%\begin{equation}\label{eq:OE}
%\lim_{m \rightarrow \infty} \lim_{\mu \rightarrow 0} \partial_\mu e_{\mu,m} = \lim_{\mu \rightarrow 0} \lim_{m \rightarrow \infty} \partial_\mu e_{\mu,m}
%\end{equation}
%as $\rho \rightarrow 0$. That is, that the formulas for the expectations of observables under soft-core potentials can be extended to that of the hard-core case as a function of the scattering length and the density.
%\end{theorem}
%The tools of the proof are mainly drawn from the strategies and the operators defined in \cite{CJL20,CJL21}.

However, for technical reasons, we will have to change the point of view slightly: instead of the setting described above in which $\rho_{\hc}$ is fixed and $e$ is computed from \eqref{eq:hsimpleq}-\eqref{eq:hsimpleq-norm}, we will fix $e$ and compute $\rho_{\hc}$ from \eqref{eq:hsimpleq}-\eqref{eq:hsimpleq-norm}.
This change of point of view was also needed in \cite{CJL20}.
The two points of view are conjectured to be equivalent \cite{CJL20}, but proving this is the case boils down to proving that $e\mapsto\rho_{\hc}(e)$ is invertible, which remains an open problem (for a fixed integrable potential, strict monotonicity is known for small and for large $e$ \cite[Theorem~1.3]{CJL21}, but the small-$e$ threshold degenerates as $m\to\infty$).
Instead, we will work at fixed $e$ and use only the following form of invertibility at small $e$, which we prove in Section~\ref{sec:LHY}: for $e\le e_0$ the density and the energy are comparable, $c\,e\le\rho_{\hc}(e)$ and $|e-2\pi\rho_{\hc}(e)|\le C e^{3/2}$, so that along the branch $e\mapsto\rho_\hc(e)$ the regimes $e\to0$ and $\rho_\hc\to0$ coincide and $e=2\pi\rho_\hc(1+o(1))$. This is what allows us, at the end of the computation, to rewrite asymptotic expansions in the small parameter $e$ as expansions in the small parameter $\rho_\hc$, which is how Theorems~\ref{thm:LHY} and~\ref{thm:cf} are stated; no inversion of $e\mapsto\rho_\hc(e)$ is needed.

%{\color{red}move?\color{blue}Throughout, $e>0$ is fixed and denotes the unperturbed energy: for every $m$, the perturbed family $\mu\mapsto e_{\mu,m}$ of \eqref{eq:sc-pert} satisfies $e_{0,m}=e$, and likewise $e_{0,\hc}=e$ for \eqref{eq:hc-pert}.}

\begin{theorem}[Condensate Fraction Expectation]\label{thm:OE}
Consider $e>0$.
\begin{enumerate}[label=(\roman*)]
\item \textup{(Energy.)} As $m\to\infty$, the densities of the soft-core converge to that of the hard-core as
\begin{equation}\label{eq:OE-i}
\bar\rho_m \longrightarrow \rho_\hc, \qquad |\bar\rho_m-\rho_\hc| = O(m^{-1/2}).
\end{equation}
Moreover, each soft-core equation satisfies the low-density expansion \eqref{eq:lhyeq}, with $a_0$ replaced by the scattering length $a_m$ of $v_m$, and the error term is uniform in $m$.
\item \textup{(Condensate fraction.)} Let $e_{\mu,m}$ and $e_{\mu,\hc}$ denote the energies of the modified equations \eqref{eq:obssimpleq}, with $v=v_m$ and $\rho=\bar\rho_m$, and \eqref{eq:obshsimpleq}. Then
\begin{equation}\label{eq:OE-ii}
\lim_{m\to\infty}\partial_\mu e_{\mu,m}\big|_{\mu=0} = \partial_\mu e_{\mu,\hc}\big|_{\mu=0},
\end{equation}
that is, the soft-core condensate fractions converge to the hard-core condensate fraction. Moreover, each soft-core equation satisfies the condensate-fraction expansion \eqref{eq:cfeq}, with $a_0$ replaced by $a_m$, and the error term can be taken uniform in $m$.
\end{enumerate}
\end{theorem}
Theorem~\ref{thm:OE} implies Theorems~\ref{thm:LHY} and~\ref{thm:cf}. For each fixed $m$, the low-density expansion of \cite{CJL20} applied to $v_m$ reads
\begin{equation}\label{eq:OE-i-exp}
\bar\rho_m \;=\; \frac{e}{2\pi a_m}\Bigl(1-\tfrac{64\sqrt2\,a_m}{15\pi}\sqrt e + o(\sqrt e)\Bigr),
\end{equation}
where $a_m = 1-\tanh(\sqrt m)/\sqrt m$ is the scattering length of $v_m$. The central point is Theorem~\ref{thm:OE} itself, proved in Theorems~\ref{thm:rate} and~\ref{thm:eta-conv}, and in particular the fact that the $o(\sqrt e)$ term may be taken \emph{uniform in} $m$: it is bounded by $\lambda(e)\sqrt e$ for some $\lambda(e)\to0$ independent of $m$, and such a bound survives the limit $m\to\infty$. Hence, by \eqref{eq:OE-i} and $a_m\to1$,
\[
\rho_\hc=\frac{e}{2\pi}\Bigl(1-\tfrac{64\sqrt2}{15\pi}\sqrt e +o(\sqrt e)\Bigr),\qquad e\to0,
\]
which inverts, via $\sqrt e=\sqrt{2\pi\rho_\hc}\,(1+o(1))$, to \eqref{eq:LHY}. Likewise the expansion of $\eta_m$ established in Section~\ref{sec:cf} is uniform in $m$, so letting $m\to\infty$ and using \eqref{eq:OE-ii} yields Theorem~\ref{thm:cf} (again via $e=2\pi\rho_\hc(1+o(1))$). This uniformity is the critical point. The errors in the fixed-$m$ expansions depend a priori on the potential (e.g.\ through $\int v_m$, which diverges as $m\to\infty$), and convergence at fixed $e$ alone does not rule out an error deteriorating as $m\to\infty$. Uniformity is obtained by arranging that $v_m$ enters the estimates only through the weight $(1-\phi_m)v_m$ --- where $\phi_m$ is the scattering solution of $v_m$, defined in Section~\ref{sec:LHY} --- whose total mass is $4\pi a_m\le4\pi$, and through the bound $|S_m|\le1$ on the structure factor \eqref{eq:Sm}. Finally, the rate in \eqref{eq:OE-i} extends uniformly in $\mu$ (Lemma~\ref{lem:uniform-energy}) and is the key quantitative input to the limit exchange \eqref{eq:OE-ii}, carried out for the condensate fraction in Theorem~\ref{thm:eta-conv}. The convergence \eqref{eq:OE-i} itself is obtained from the $L^1$ convergence $\int|\sigma_m-\tau|\to0$ via $|\bar\rho_m-\rho_\hc|=\rho_\hc\,\bar\rho_m\bigl|\int(\sigma_m-\tau)\bigr|$.

\begin{remark}[Other observables]\label{rem:linear}
We have restricted the discussion to the condensate fraction, but the framework extends immediately to other observables, in particular to the momentum distribution. For any one- or two-particle observable $\mathcal U$, the Hellmann--Feynman procedure of Section~\ref{sec:observables}, with $\frac1N\sum_jP_j$ replaced by $\mathcal U$, yields a modified Simple Equation \cite{CJL21} in which the perturbation enters \emph{affinely} in the solution, as $\mu(\beta+\mathfrak{L}\,u_\mu)$ for a source $\beta$ and a linear operator $\mathfrak{L}$ determined by $\mathcal U$ and independent of $\mu$: this is because the expectation is a first-order variation of the ground-state energy, and the Simplified Approach closes the correlation hierarchy at the two-point function \cite{Li63,J25}. For the condensate fraction, $\beta=0$ and $\mathfrak{L}=-2\,\mathrm{Id}$, i.e.\ the shift $-\Delta\mapsto-\Delta+2\mu$ in \eqref{eq:obssimpleq}; for the momentum distribution, $\mathfrak{L}=0$ and $\beta$ is proportional to $\hat u_0(k)\cos(k\cdot x)$ \cite[(1.20)]{CJL21}. The limit-exchange argument of Section~\ref{sec:cf} depends on the observable only through the verification of Lemmas~\ref{lem:uniform-energy}--\ref{lem:hc-diff} for the perturbed family, so extending part (ii) of Theorem~\ref{thm:OE} to the momentum distribution, or to any other such observable, reduces to verifying these three lemmas for its modified equation.
\end{remark}

We will now discuss the strategy to prove Theorem~\ref{thm:OE}. The Simple Equation, hard-core or not, is extremely difficult to tackle directly due to the convolutional non-linearity. For a given potential, $v(x)$, the strategy that was leveraged previously \cite{CJL20,CJL21} was to use the following iteration scheme
(we will introduce other iteration schemes, and to differentiate them, we will use different letters for the iterates; the soft-core iterates from \cite{CJL20,CJL21} will be denoted by $u_n$)
\begin{align}
	(-\Delta +4e +v(|x|)) u_n(|x|) & = v(|x|) + 2 e \rho_{n-1} u_{n-1} \ast u_{n-1} (|x|) \\ 
	\rho_n & =  \frac{2e}{ \int (1-u_n)v}
\end{align}
with $u_0=0$. $u_n$ converges to the solution $u$ of the Simple Equation in Eq.~\eqref{eq:simpleq} for potential $v$, see \cite[Theorem 1.3]{CJL20}.
%Nice monotonic behaviors were proven that $u_{n-1} \le u_n$ and $\rho_{n-1} \le \rho_n$. These properties make it fairly easy to apply harmonic analysis within an inductive construction since $\rho_{n-1} u_{n-1} \ast u_{n-1} \le \rho_n u_n \ast u_n$. 
This iteration, for the soft-core potential $v_m\equiv m\mathbbm1_{\Bone}$ is:
\begin{equation}\label{eq:softcore}
	(-\Delta + 4e + v_m)\sigma_{m,n} = v_m + 2e\bar{\rho}_{m,n-1}\sigma_{m,n-1} * \sigma_{m,n-1}, \qquad
	\bar{\rho}_{m,n} = \frac{2e}{m\int_{\Bone}(1 - \sigma_{m,n})}.
\end{equation}
with $\sigma_{m,0} = 0$.
%$m$ and $n$ are two independent indices and it is evident that many of the results for a fixed $m$ can be drawn from \cite{CJL20,CJL21}.
By \cite[Theorem 1.3]{CJL20}, $\lim_{n \rightarrow \infty} \sigma_{m,n} \equiv \sigma_m$ is the solution to Eq.~\eqref{eq:ssimpleq}; the properties of these iterates that we use are collected in Lemma~\ref{lem:softcore-props}.

We will use the following iteration for the hard-core equation which follows from \eqref{eq:hsimpleq} and \eqref{eq:ehard}:
\begin{align}
\Delta \tau_n & = 4e\tau_n - 2e\rho_{n-1}\tau_{n-1} * \tau_{n-1}, \quad &&|x| > 1, \label{eq:iterate}\\
\tau_n &= 1, \quad &&|x| \le 1,
\end{align}
with
\begin{equation}\label{eq:rhon}
\rho_n = \frac{2e}{2e\bigl[\frac{8\pi}3 - \rho_{n-1}\int_{\Bone} \tau_{n-1} * \tau_{n-1}\bigr] - 4\pi\,\partial_r \tau_n\big|_{1^+}}.
\end{equation}
We initialize the iteration with $\Delta \tau_1 = 4e\tau_1$ on $|x|>1$, $\tau_1=1$ on $\Bone$, that is,
\begin{equation}\label{eq:u1}
\tau_1(r) = \begin{cases} 1 & r \leq 1 \\ e^{-2\sqrt{e}(r-1)}/r & r > 1 \end{cases}, \qquad
\rho_1 = \frac{1}{2 \int \tau_1} = \frac{2e}{4\pi(1 + 2\sqrt{e}) + \frac{16\pi}3e}.
\end{equation}
To show that $\tau_n$ converges to the solution of \eqref{eq:hsimpleq}, we will prove that $\tau_{n}$ and $\rho_n$ are increasing in $n$ and use the monotone convergence theorem, see Lemmas~\ref{lem:properties} and~\ref{lem:convergence}; uniqueness of the limit is Lemma~\ref{lem:uniqueness}. The monotonicity and comparison arguments throughout rest on the maximum principle in the form of Lemma~\ref{lem:maxp} and Corollary~\ref{cor:maxp}.

%Whenever we refer to $\tau_n$ or $\lim_{n \rightarrow \infty} \tau_n \equiv \tau$, we are henceforth refering to the hard-core Simple Equation.

\begin{figure}[htbp]
\centering
\includegraphics[width=\textwidth]{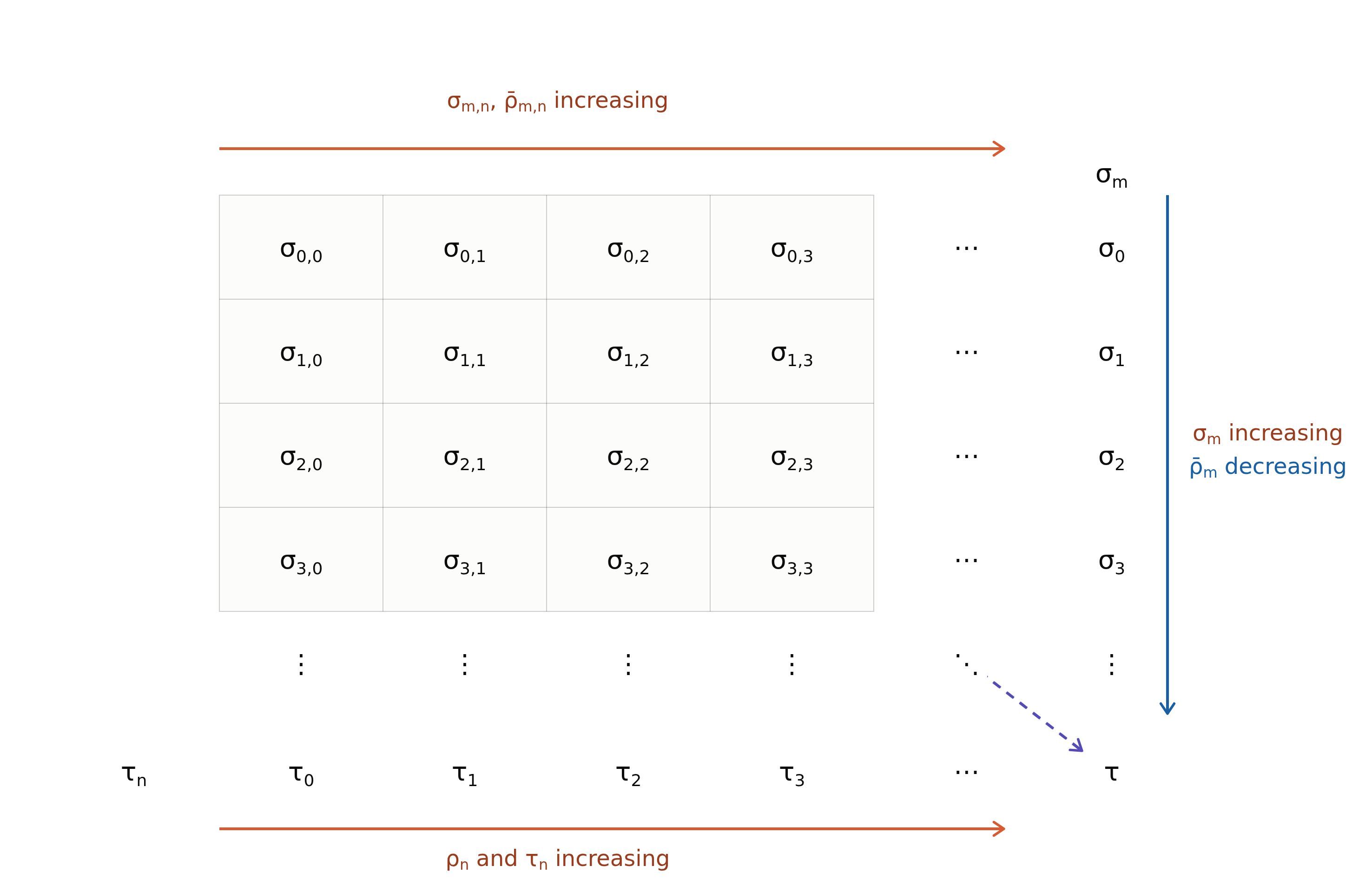}
\caption{ The soft-core iterates $\sigma_{m,n}$, arranged by the height $m$ of the potential and the iteration index $n$. Increasing $m$ raises the iterates inside the core while lowering the densities $\bar\rho_{m,n}$, which obstructs a direct comparison along the diagonal $\sigma_{n,n}$.}
\label{fig:sigma}
\end{figure}

The main idea is then to compare $\tau_n$ and $\sigma_{m,n}$ (for each fixed $n$, Lemma~\ref{lem:per-n} shows that $\sigma_{m,n}\to\tau_n$ as $m\to\infty$).
However, we cannot easily do this directly, as we will show that, keeping $e$ fixed, as one increases $m$, $\bar{\rho}_{m}$ decreases (Lemma~\ref{lem:rhom-decreasing}).
Simultaneously an increase in the potential in $\Bone$ means that $\sigma_{m}$ must increase with $m$ and converge to 1 as $m\to\infty$. This makes comparison of terms $\bar{\rho}_{m,n} \sigma_{m,n} \ast \sigma_{m,n}$ and $\bar{\rho}_{m-1,n} \sigma_{m-1,n} \ast \sigma_{m-1,n}$ difficult since the decrease in the density must be balanced with any potential increase in the solution with iteration. Naturally this difficulty carries over to the iterates $\sigma_{m,n}$ and $\bar{\rho}_{m,n}$. In addition, this completely impedes the comparison of terms along the diagonal of Figure~\ref{fig:sigma}.

In order to work around this impasse, we will avoid the competition between inverse monotonicities, by introducing two more iterations.

The type~$\gamma$ iterates $\gamma_n$ are defined by $\gamma_0 = 0$ and, for $n\ge 1$,
\begin{equation}\label{eq:typeI}
(-\Delta + 4e + v_n)\gamma_n = v_n + 2e\rho_{n-1}\gamma_{n-1} * \gamma_{n-1},
\end{equation}
where $v_n = n \mathbbm{1}_{\Bone}$ and the $\rho_{n-1}$ are from the hard-core iteration \eqref{eq:rhon}.
In other words, we use the sequence of $\rho_n$'s from the hard-core iteration, and increase the potential at the same time as we proceed through the iteration.

\begin{figure}[htbp]
\centering
\includegraphics[width=\textwidth]{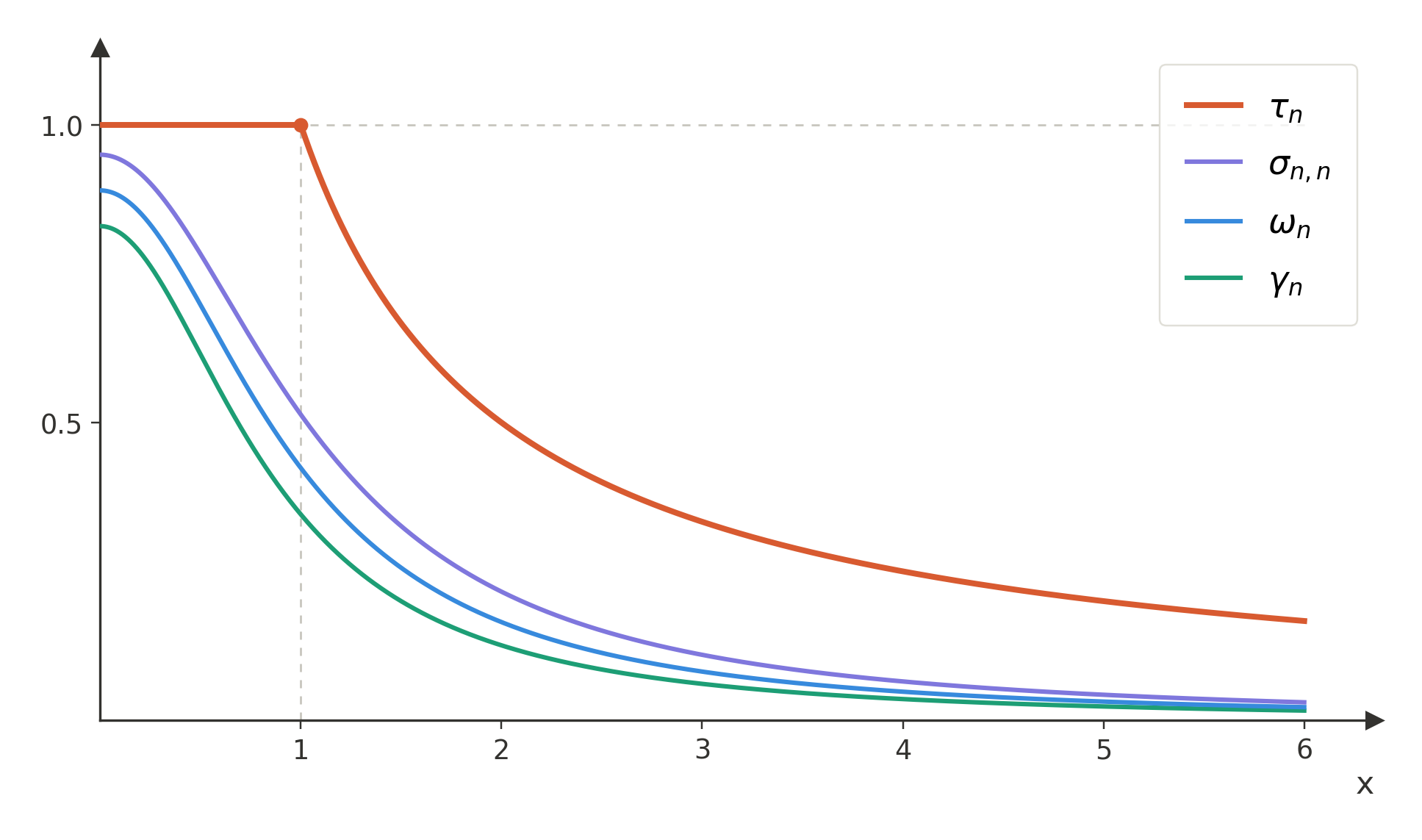}
\caption{ The nesting of the iterates: $\gamma_n\le\omega_n\le\tau$ and $\gamma_n\le\tau_n\le\tau$, with $\omega_n\le\sigma_{n,n}\le\sigma_n$; all of them converge to the hard-core solution $\tau$ in $L^1$.}
\label{fig:nested}
\end{figure}

The type~$\omega$ iterates $\omega_n$ are defined by $\omega_0 = 0$ and
\begin{equation}\label{eq:typeII}
(-\Delta + 4e + v_n)\omega_n = v_n + 2e\rho_{\hc}\omega_{n-1} * \omega_{n-1}.
\end{equation}
where $\rho_{\hc}$ is the hard-core density. 
The difference with the $\gamma$ iteration is only that $\rho_{n-1}$ is replaced by its limit $\rho_{\hc}$. The relations between the various iterates are summarized in Figure~\ref{fig:nested}.

Type $\gamma$ iterates will allow us to directly compare $\rho_{n} \tau_n \ast \tau_n$ to $\rho_n \gamma_n \ast \gamma_n$ since the densities are the same and we will be able to prove that $\tau_n \ge \gamma_n$ (Lemma~\ref{lem:gamma-bounds}) and $\gamma_n \rightarrow \tau$ in $L^1$ (Lemma~\ref{lem:gamma-limit}). This will allow us to create an envelope such that $\gamma_n \le \omega_n \le \tau$ (Lemma~\ref{lem:omega}). We cannot easily compare $\omega_n$ to $\tau_n$ since $\rho_\hc \ge \rho_{n} \; \forall n$. The lower bound will be obtained by comparing the non-linear terms $\rho_{\hc} \omega_n \ast \omega_n$ and $\rho_{n} \gamma_n \ast \gamma_n$ since $\rho_{\hc} \ge \rho_{n}$. As a result we will show that $\omega_n \overset{L^1}{\rightarrow} \tau$ (Corollary~\ref{cor:omega-u}). Ultimately we will then lower bound $\sigma_{n,n}$ with $\omega_n$ (Lemma~\ref{lem:comparison}) which will enable us to show that $\sigma_n \overset{L^1}{\rightarrow} \tau$ since $\sigma_{n,n} \le \sigma_n$ (Theorem~\ref{thm:convergencesigma}). This strategy will provide $L^1$ convergence which is necessary but insufficient since the rate of convergence is also critical to verifying the LHY formula for the hard-core Simple Equation. To do this will require carefully studying the rate of $\bar{\rho}_m$ to $\rho_\hc$, which will be analytically estimated in Theorem~\ref{thm:rate}, using the derivative formula \eqref{eq:rhoderiv} of Lemma~\ref{lem:rhom-decreasing}. The LHY formula, Theorem~\ref{thm:LHY}, is then proved in Section~\ref{sec:LHY} from the identity of Lemma~\ref{lem:e-minus-rho} and the uniform bounds of Lemmas~\ref{lem:Sm-bound} and~\ref{lem:eta-w-bound}; the condensate fraction, Theorem~\ref{thm:cf}, is proved in Section~\ref{sec:cf} through Lemmas~\ref{lem:ND}--\ref{lem:hc-diff}, Theorem~\ref{thm:eta-conv} and Corollary~\ref{cor:cf-bridge}.

%=====================================================================
\section{Iterative Construction of the Hard-Core Solution}
%=====================================================================

%Define the sequence which for $n \geq 1$
%\begin{align}
%\Delta u_n & = 4eu_n - 2e\rho_{n-1}u_{n-1} * u_{n-1}, \quad &&|x| > 1, \label{eq:iterate}\\
%u_n &= 1, \quad &&|x| \le 1,
%\end{align}
%with
%\begin{equation}\label{eq:rhon}
%\rho_n = \frac{2e}{2e\bigl[2\Bvol - \rho_{n-1}\int_{\Bone} u_{n-1} * u_{n-1}\bigr] - |S_1(0)|\,\partial_r u_n\big|_{1^+}}.
%\end{equation}
%The first iterate satisfies $\Delta u_1 = 4eu_1$ on $|x|>1$, $u_1=1$ on $\Bone$, giving
%\begin{equation}\label{eq:u1}
%u_1(r) = \begin{cases} 1 & r \leq 1 \\ e^{-2\sqrt{e}(r-1)}/r & r > 1 \end{cases}, \qquad
%\rho_1 = \frac{2e}{|S_1(0)|(1 + 2\sqrt{e}) + 4e\Bvol}.
%\end{equation}

Before constructing the iterative solution scheme to the hard-core, we will introduce two lemmas will be used repeatedly in the following.
%The arguments below compare a function against a barrier on the set where the desired inequality fails.

\begin{lemma}[Comparison principle]\label{lem:maxp}
Let $D\subseteq\R^3$ be open and let $0\le W\in L^\infty(D)$. Let $w$ be continuous on $\overline D$ and satisfy 
\begin{equation}\label{eq:maxp}
(-\Delta+W)\,w \ \ge\ 0 \quad\text{on }D
\end{equation}
in the sense of distributions. If $w\ge0$ on $\partial D\equiv \overline D\setminus D$ and $\liminf_{|x|\to\infty}w(x)\ge0$ (when $D$ is unbounded), then $w\ge0$ on $D$.
\end{lemma}

\begin{proof}
Fix $\varepsilon>0$ and let $\Sigma_\varepsilon=\{x\in D: w(x)<-\varepsilon\}$, an open set. Since $\liminf_{|x|\to\infty}w(x)\ge0$, so $\Sigma_\varepsilon$ is bounded. Since $w\ge0>-\varepsilon$ on $\partial D$, we have $\overline{\Sigma_\varepsilon}\cap\partial D=\emptyset$, hence $\overline{\Sigma_\epsilon}$ is a compact subset of $D$, and $w=-\varepsilon$ on $\partial\Sigma_\varepsilon$. Suppose $\Sigma_\varepsilon\neq\emptyset$. On $\Sigma_\varepsilon$ we have $w<0$ and $W\ge0$, so \eqref{eq:maxp} gives $\Delta w\le W w\le0$: $w$ is superharmonic on $\Sigma_\varepsilon$,
%(mollification preserves $\Delta w\le0$, and the mean-value inequality passes to the locally uniform limit),
so it attains its minimum on $\partial\Sigma_\varepsilon$, where $w=-\varepsilon$, contradicting the definition of $\Sigma_\varepsilon$. Hence $\Sigma_\varepsilon=\emptyset$ for every $\varepsilon>0$.
\end{proof}

\begin{corollary}[Constant barriers]\label{cor:maxp}
Let $D\subset\R^3$ be open and bounded, $c>0$, $0\le W\in L^\infty(D)$, and let $w$ be continuous on $\overline D$ with
\[
0\ \le\ (-\Delta+c+W)\,w\ \le\ F \quad\text{on }D
\]
in the sense of distributions, for a constant $F\ge0$. Then, on $D$,
\[
\min\Bigl\{0,\ \min_{\partial D}w\Bigr\}\ \le\ w\ \le\ \max\Bigl\{0,\ \max_{\partial D}w,\ \tfrac{F}{c}\Bigr\}.
\]
\end{corollary}

\begin{proof}
Set $\alpha=\min\{0,\min_{\partial D}w\}$ and $M=\max\{0,\max_{\partial D}w,F/c\}$. We apply Lemma~\ref{lem:maxp} to $w-\alpha$: since $\alpha\le0$, $(-\Delta+c+W)(w-\alpha)\ge-(c+W)\alpha\ge0$, and $w-\alpha\ge0$ on $\partial D$.
Therefore, $w \ge \alpha$.
Next, we apply Lemma~\ref{lem:maxp} to $M-w$: since $M\ge F/c$ and $W\ge0$, $(-\Delta+c+W)(M-w)\ge cM-F\ge0$, and $M-w\ge0$ on $\partial D$.
Thus, $w \le M$.
\end{proof}

With this lemma and corollary proven, we can now prove several simple but useful properties of the iterative solutions to the hard-core Simple Equation. For $n \geq 2$, the solution is obtained via the Green's function of $(-\Delta + 4e)$ on $\mathbb{R}^3$, namely $Y_{4e}(x) = e^{-2\sqrt{e}|x|}/(4\pi|x|)$. Solving Eq.~\eqref{eq:iterate} directly, we find
\begin{equation}\label{eq:un-exterior}
\tau_n(r) = a_n \frac{e^{-2\sqrt{e}(r-1)}}{r} + \frac{\sqrt{e}\,\rho_{n-1}}{2r}\int_1^\infty s[\tau_{n-1} * \tau_{n-1}](s)\,e^{-2\sqrt{e}|s-r|}\,ds, \qquad r > 1,
\end{equation}
with
\begin{equation}\label{eq:an}
a_n = 1 - \frac{\sqrt{e}\,\rho_{n-1}}{2}\int_1^\infty s[\tau_{n-1} * \tau_{n-1}](s)\,e^{-2\sqrt{e}(s-1)}\,ds.
\end{equation}
Therefore each $\tau_n$ is continuous on $\mathbb{R}^3$ and smooth on $\Bone$ and $B_1^c$ separately. Differentiating \eqref{eq:un-exterior} at $r = 1^+$ gives
\begin{equation}\label{eq:dun}
\partial_r \tau_n\big|_{1^+} = -1 - 2\sqrt{e} - 4\sqrt{e}\,a_n.
\end{equation}

\begin{lemma}[Properties of hard-core iterates]\label{lem:properties}
For all $n \geq 1$:
\begin{enumerate}[label=(\roman*)]
\item $0 \leq \tau_n \leq 1$,
\item $\int \tau_n \leq 1/\rho_n$,
\item $\tau_{n} \leq \tau_{n+1}$ and $\rho_{n} \leq \rho_{n+1}$.
\end{enumerate}
\end{lemma}

\begin{proof}
We proceed by induction. For $n=1$, (i) and (ii) are clear from \eqref{eq:u1}. To see $\tau_1\le \tau_2$, we write the difference using the explicit formula and observe the integrand is non-negative by the triangle inequality. The inequality $\rho_2\ge\rho_1$ then follows from the comparison of the denominators in \eqref{eq:rhon}.

Assume the statements hold for all indices up to $n-1$.
\newline
%\textit{Non-negativity.} If $\Sigma^- = \{x\in B_1^c : \tau_n(x)<0\}\neq\emptyset$, then on $\Sigma^-$ we have $\Delta \tau_n = 4e\tau_n - 2e\rho_{n-1}\tau_{n-1}*\tau_{n-1} < 0$, so $\tau_n$ is superharmonic and attains its minimum on $\partial\Sigma^-$, where $\tau_n=0$, contradiction.
%\textit{Non-negativity.} Apply Lemma~\ref{lem:maxp} with $\Sigma^- = \{x\in B_1^c : \tau_n(x)<0\}\neq\emptyset$. Since on $\Sigma^-$ we have $\Delta \tau_n = 4e\tau_n - 2e\rho_{n-1}\tau_{n-1}*\tau_{n-1} < 0$ it follows that $\tau_n  0$.
\textit{\textbf{Non-negativity.}} On $\{|x|>1\}$ we have $(-\Delta+4e)\tau_n = 2e\rho_{n-1}\,\tau_{n-1}*\tau_{n-1} \ge 0$ by the induction hypothesis, while $\tau_n = 1$ on $\partial B_1$ and $\tau_n\to0$ at infinity by \eqref{eq:un-exterior}. Lemma~\ref{lem:maxp} with $D=\{|x|>1\}$, $c=4e$, $W=0$ gives $\tau_n\ge0$.
\newline
\textit{\textbf{Upper bound} $\tau_n\le1$.} Since $(-\Delta+4e)(1-\tau_n) = 4e - 2e\rho_{n-1}\,\tau_{n-1}*\tau_{n-1} \ge 2e > 0$ by $\tau_{n-1}\le1$ and $\rho_{n-1}\int\tau_{n-1}\le1$, while $1-\tau_n = 0$ on $\partial B_1$ and $1-\tau_n\to1$ at infinity, Lemma~\ref{lem:maxp} gives $\tau_n\le1$.
%\textit{Upper bound $\tau_n\le 1$.} Suppose $\Sigma^+ = \{x\in B_1^c : \tau_n(x)>1\}\neq\emptyset$. On $\Sigma^+$,
%\[
%\Delta \tau_n = 4e\tau_n - 2e\rho_{n-1}\tau_{n-1}*\tau_{n-1} > 4e - 2e\rho_{n-1}\int \tau_{n-1} \ge 4e - 2e = 2e > 0,
%\]
%using the induction hypothesis $\rho_{n-1}\int \tau_{n-1}\le 1$. Hence $\tau_n$ is subharmonic on $\Sigma^+$ and attains its maximum on $\partial\Sigma^+$, where $\tau_n=1$, contradiction.
\newline
\textit{\textbf{Monotonicity.}} For $\tau_n\le \tau_{n+1}$, observe that on $|x|>1$,
\[
(-\Delta+4e)(\tau_{n+1}-\tau_n) = 2e[\rho_n \tau_n*\tau_n - \rho_{n-1}\tau_{n-1}*\tau_{n-1}] \ge 0,
\]
by the induction hypothesis $\rho_n\ge\rho_{n-1}$ and $\tau_n\ge \tau_{n-1}$. Since $\tau_{n+1}-\tau_n=0$ on $\partial B_1$, the maximum principle yields $\tau_{n+1}\ge \tau_n$. This implies $\partial_r \tau_{n+1}|_{1^+} \ge \partial_r \tau_n|_{1^+}$, and together with the induction hypothesis the denominator in \eqref{eq:rhon} does not increase, so $\rho_{n+1}\ge\rho_n$.
\newline
\textit{\textbf{Integral bound.}} Integrating \eqref{eq:iterate} over $\{|x|>1\}$ and using the definition of $\rho_n$ gives
\begin{equation}\label{eq:balance}
2\int \tau_n = \rho_{n-1}\Bigl(\int \tau_{n-1}\Bigr)^2 + \frac{1}{\rho_n} \le \int \tau_{n-1} + \frac{1}{\rho_n} \le \int \tau_n + \frac{1}{\rho_n},
\end{equation}
hence $\int \tau_n \le 1/\rho_n$.
\end{proof}

\begin{lemma}[Convergence]\label{lem:convergence}
The limits $\tau = \lim_{n \to \infty} \tau_n$ and $\rho_{\hc} = \lim_{n \to \infty} \rho_n$ exist, satisfy $0 \leq \tau \leq 1$, $\rho_{\hc}\int \tau = 1$, and $\|\tau_n - \tau\|_{L^1} \to 0$.
\end{lemma}

\begin{proof}
$\tau_n$ is increasing and bounded by 1, so $\tau_n \to \tau$ pointwise. Since $0 \le \tau_n \le \tau$ and $\int \tau_n \le 1/\rho_1 < \infty$, monotone convergence gives $\tau \in L^1$ and $\|\tau - \tau_n\|_{L^1} = \int (\tau - \tau_n) \to 0$. The density $\rho_n$ is increasing and bounded (by $\rho_n \le 1/\int \tau_1$), so $\rho_n \to \rho_{\hc}$. Passing the integral relation $\rho_n\int \tau_n = 1$ (which follows from \eqref{eq:balance} and the monotone limits) to the limit gives $\rho_{\hc}\int \tau = 1$. The fact that $\tau$ solves \eqref{eq:hsimpleq} follows by taking $n\to\infty$ in the weak formulation of \eqref{eq:iterate}.
\end{proof}

\begin{lemma}[Uniqueness]\label{lem:uniqueness}
$(\tau, \rho_{\hc})$ is the unique solution of \eqref{eq:hsimpleq} with $0 \leq \tau \leq 1$ and $\tau \in L^1$.
\end{lemma}

\begin{proof}
Let $\tau'$ be another such solution with density $\rho'$. We show $\tau' \geq \tau_n$ for all $n$ by induction. Assume $\tau' \geq \tau_{n-1}$. If $\Upsilon = \{|x| > 1 : \tau'(x) < \tau_n(x)\} \neq \emptyset$, then on $\Upsilon$,
\[
\Delta(\tau' - \tau_n) = 4e(\tau' - \tau_n) - 2e[\rho'\tau' * \tau' - \rho_{n-1}\tau_{n-1} * \tau_{n-1}] < 0,
\]
so $\tau' - \tau_n$ is superharmonic on $\Upsilon$ and attains its minimum on $\partial\Upsilon$, where it vanishes (since $\tau'=\tau_n=1$ on $\partial B_1$). Contradiction. Therefore $\tau' \geq \tau_n$ for all $n$, hence $\tau' \geq \tau$ and $\rho' \geq \rho_{\hc}$. But $\int \tau' = 1/\rho'$ and $\int \tau = 1/\rho_{\hc}$, so $\rho' \leq \rho_{\hc}$. Thus $\rho' = \rho_{\hc}$ and $\tau' = \tau$ a.e.
\end{proof}

%=====================================================================
\section{Intermediate Iterative Systems}\label{sec:intermediate}
%=====================================================================

%The type~$\gamma$ iterates $\gamma_n$ are defined by $\gamma_0 = 0$ and, for $n\ge 1$,
%\begin{equation}\label{eq:typeI}
%(-\Delta + 4e + v_n)\gamma_n = v_n + 2e\rho_{n-1}\gamma_{n-1} * \gamma_{n-1},
%\end{equation}
%where $v_n = n \mathbbm{1}_{\Bone}$ and the $\rho_{n-1}$ are from the hard-core iteration.
%
%The type~$\omega$ iterates $\omega_n$ are defined by $\omega_0 = 0$ and
%\begin{equation}\label{eq:typeII}
%(-\Delta + 4e + v_n)\omega_n = v_n + 2e\rho_{\hc}\omega_{n-1} * \omega_{n-1}.
%\end{equation}
\subsection{The $\gamma_n$ iteration}
We recall the definition of the $\gamma_n$ iteration \eqref{eq:typeI}:
\begin{equation}\label{eq:gamma2}
\gamma_0=0
,\quad
(-\Delta + 4e + v_n)\gamma_n = v_n + 2e\rho_{n-1}\gamma_{n-1} * \gamma_{n-1}.
\end{equation}
In this section we will prove that $\gamma_n$ converges to $\tau$ in $L^1$ norm.
\begin{lemma}\label{lem:gamma-bounds}
For all $n \geq 0$,
\begin{enumerate}[label=(\roman*)]
  \item $0 \leq \gamma_n \leq 1$,
  \item $\gamma_n \leq \tau_n$,
  \item $\gamma_n \leq \gamma_{n+1}$.
\end{enumerate}
\end{lemma}

\begin{proof}
Take $\gamma_0=\tau_0=0$. By induction; the case $n=0$ is trivial. Assume $0\le\gamma_{n-1}\le1$ and $\gamma_{n-1}\le\tau_{n-1}$.
\newline
\textit{\textbf{Non-negativity.}} The right-hand side of \eqref{eq:typeI} is non-negative and $(-\Delta+4e+v_n)^{-1}$ is positivity preserving, so $\gamma_n\ge0$.
\newline
\textit{\textbf{Upper bound $\gamma_n\le1$.}} On $\R^3$,
\[
(-\Delta+4e+v_n)(1-\gamma_n) = 4e - 2e\rho_{n-1}\gamma_{n-1}*\gamma_{n-1} \ge 4e-2e = 2e > 0,
\]
using $\gamma_{n-1}\le1$ and $\rho_{n-1}\int\gamma_{n-1}\le\rho_{n-1}\int\tau_{n-1}\le1$ (which follows from the induction hypothesis and Lemma~\ref{lem:properties}); moreover $1-\gamma_n\to1$ at infinity. Lemma~\ref{lem:maxp} with $D=\R^3$, $c=4e$, $W=v_n$ gives $\gamma_n\le1$.
\newline
\textit{\textbf{Monotonicity.}} Assuming $\gamma_{n-1}\le\gamma_n$, on $\R^3$,
\begin{align*}
&(-\Delta+4e+v_{n+1})(\gamma_{n+1}-\gamma_n) \\
&\quad = (v_{n+1}-v_n)(1-\gamma_n) + 2e[\rho_n\gamma_n*\gamma_n - \rho_{n-1}\gamma_{n-1}*\gamma_{n-1}] \ge 0,
\end{align*}
since $v_{n+1}\ge v_n$, $\gamma_n\le1$ (by the upper bound above), $\rho_n\ge\rho_{n-1}$ (Lemma~\ref{lem:properties}), and $\gamma_n\ge\gamma_{n-1}\ge0$; moreover $\gamma_{n+1}-\gamma_n\to0$ at infinity. Lemma~\ref{lem:maxp} with $D=\R^3$, $c=4e$, $W=v_{n+1}$ gives $\gamma_{n+1}\ge\gamma_n$.
\newline
\textit{\textbf{Comparison $\gamma_n\le\tau_n$.}} On $\Bone$ this reads $\gamma_n\le1=\tau_n$. On $\{|x|>1\}$,
\[
(-\Delta+4e)(\tau_n-\gamma_n) = 2e\rho_{n-1}[\tau_{n-1}*\tau_{n-1}-\gamma_{n-1}*\gamma_{n-1}] \ge 0
\]
by the induction hypothesis, while $\tau_n-\gamma_n = 1-\gamma_n(1)\ge0$ on $\partial B_1$ and $\tau_n-\gamma_n\to0$ at infinity. Lemma~\ref{lem:maxp} with $D=\{|x|>1\}$, $c=4e$, $W=0$ gives $\gamma_n\le\tau_n$.

%We first prove (i) and (ii) by induction. Take $\gamma_0=\tau_0=0$.
%Assume $0\le\gamma_{n-1}\le 1$ and $\gamma_{n-1}\le\tau_{n-1}$.
%By lemma \ref{lem:properties}, the right-hand side of \eqref{eq:gamma2} is non-negative, and $(-\Delta+4e+v_n)^{-1}$ is positivity preserving, so $\gamma_n\ge0$. For the upper bound, suppose that $\Sigma^+ := \{x : \gamma_n(x) > 1\}\neq\emptyset$. On $\Sigma^+$,
%\[
%\Delta\gamma_n = 4e\gamma_n - 2e\rho_{n-1}\gamma_{n-1}*\gamma_{n-1} +v_n(\gamma_n-1) > 4e - 2e\rho_{n-1}\int\gamma_{n-1}
%\]
%Furthermore, $\int \gamma_{n-1} \le \int\tau_{n-1}$, and so, by lemma \ref{lem:properties}, $\int \gamma_{n-1}\le1/\rho_{n-1}$, and so
%\[
%\Delta\gamma_n
%\ge 2e > 0,
%\]
%so $\gamma_n$ is subharmonic and attains its maximum on $\partial\Sigma^+$, where $\gamma_n= 1$, which is a contradiction.
%Therefore $0\le \gamma_n \le 1$.
%To show that $\gamma_n \le \tau_n$, we write, for $|x|>1$,
%\[
%(-\Delta+4e)(\tau_n-\gamma_n) = 2e\rho_{n-1}[\tau_{n-1}*\tau_{n-1} - \gamma_{n-1}*\gamma_{n-1}] \ge 0,
%\]
%and for $|x|=1$, $\tau_n-\gamma_n = 1-\gamma_n(1) \ge 0$.
%Proceeding as above {\color{red}details}, using the Maximum principle, we find gives $\tau_n\ge\gamma_n$.
%
%Finally, we turn to (iii), which we also prove by induction.
%The base case is trivial. Assume $\gamma_{n-1} \le \gamma_{n}$.
%We have
%\begin{align*}
%&(-\Delta+4e+v_{n+1})(\gamma_{n+1}-\gamma_n) \\
%&\quad = (v_{n+1}-v_n)(1-\gamma_n) + 2e[\rho_n\gamma_n*\gamma_n - \rho_{n-1}\gamma_{n-1}*\gamma_{n-1}] \ge 0,
%\end{align*}
%hence $\gamma_{n+1}\ge\gamma_n$.
\end{proof}

\begin{lemma}\label{lem:gamma-limit}
$\gamma_n \uparrow \tau$ pointwise and in $L^1$.
\end{lemma}

\begin{proof}
By Lemma~\ref{lem:gamma-bounds}, the sequence $(\gamma_n)$ is non-decreasing and bounded above by $\tau$, so $\gamma_n\uparrow\gamma\le\tau$ pointwise. Since $\tau\in L^1$ (Lemma~\ref{lem:convergence}), monotone convergence gives $\gamma\in L^1$ and
\[
\|\gamma-\gamma_n\|_{L^1} = \int(\gamma-\gamma_n) \to 0.
\]
The rest of the proof consists in showing that $\gamma=\tau$.

\textit{\textbf{Boundary value.}} Integrating \eqref{eq:gamma2} over $\mathbb{R}^3$ yields
\begin{equation}\label{eq:gamma-int}
n\int_{\Bone}(1-\gamma_n) = 4e\int\gamma_n - 2e\rho_{n-1}\|\gamma_{n-1}\|_{L^1}^2 \le 4e\int\tau_n \le \frac{4e}{\rho_1},
\end{equation}
using $\gamma_n\le\tau_n$ (Lemma~\ref{lem:gamma-bounds}) and $\int\tau_n\le1/\rho_n\le1/\rho_1$ (Lemma~\ref{lem:properties}); hence $\int_{\Bone}(1-\gamma_n) \le C/n \to 0$. To see that $\gamma_n(1)\to1$, set $\iota_n = 1-\gamma_n \ge 0$ on $\Bone$; it satisfies, for $x\in\Bone$,
\[
(-\Delta+\kappa^2)\iota_n = 4e-2e\rho_{n-1}(\gamma_{n-1}\ast\gamma_{n-1}) \equiv f_n,
\qquad \kappa=\sqrt{4e+n},
\]
with $0\le f_n\le 4e$: indeed $2e\rho_{n-1}(\gamma_{n-1}\ast\gamma_{n-1}) \le 2e\rho_{n-1}\|\gamma_{n-1}\|_{L^\infty}\|\gamma_{n-1}\|_{L^1} \le 2e\rho_{n-1}\int\tau_{n-1}\le 2e$ by Lemmas~\ref{lem:gamma-bounds} and~\ref{lem:properties}. Decompose $\iota_n = \tilde\iota_n + h_n$, where $\tilde\iota_n|_{\partial B_1}=0$ and $h_n$ is the homogeneous solution: $(-\Delta+\kappa^2)h_n=0$ with $h_n(1)=\iota_n(1)$, i.e.\ $h_n(r) = \iota_n(1)\frac{\sinh(\kappa r)}{r\sinh\kappa}$. Since $(-\Delta+\kappa^2)h_n=0$, the function $\tilde\iota_n$ satisfies the same equation as $\iota_n$, with $\tilde\iota_n=0$ on $\partial B_1$; Corollary~\ref{cor:maxp} with $D=\Bone$, $c=\kappa^2$, $W=0$, $F=4e$ gives $0\le\tilde\iota_n\le4e/\kappa^2$. Furthermore,
\[
\int_{\Bone}h_n = \frac{4\pi\iota_n(1)}{\kappa^2\sinh\kappa}(\kappa\cosh\kappa-\sinh\kappa) \sim \frac{4\pi\iota_n(1)}{\kappa}\quad(\kappa\to\infty).
\]
Since $\iota_n = \tilde\iota_n + h_n \ge h_n \ge 0$ and $\int_{\Bone}\iota_n = O(1/n)$ by \eqref{eq:gamma-int}, we obtain $\iota_n(1) = O(\kappa/n) = O(1/\sqrt{n})$, hence $\gamma_n(1)\to1$. Since $\gamma$ is the pointwise limit, $\gamma(1) = \lim_n\gamma_n(1) = 1$; and by \eqref{eq:gamma-int}, $\int_{\Bone}(1-\gamma) = \lim_n\int_{\Bone}(1-\gamma_n) = 0$, so $\gamma=1$ a.e.\ on $\Bone$.
\newline

\textit{\textbf{Continuity at the boundary.}}
$\gamma$ is continuous up to $\partial B_1$: for every $n$ and $r>1$, $\gamma(r)\ge\gamma_n(r)$, so by continuity of $\gamma_n$,
\[
\liminf_{r\to1^+}\gamma(r) \;\ge\; \gamma_n(1) \;=\; 1-O(n^{-1/2}) \quad\text{for every } n,
\]
whence $\liminf_{r\to1^+}\gamma(r)\ge1$; combined with $\gamma\le1$ this gives $\gamma(r)\to1=\gamma(1)$.
\newline

\textit{\textbf{Limit equation.}} On $\{|x|>1\}$, $v_n=0$, so \eqref{eq:gamma2} reads $(-\Delta+4e)\gamma_n = 2e\rho_{n-1}(\gamma_{n-1}\ast\gamma_{n-1})$. The convolutions converge uniformly:
\[
\|\gamma_{n-1}\ast\gamma_{n-1} - \gamma\ast\gamma\|_{L^\infty} \le (\|\gamma_{n-1}\|_{L^\infty}+\|\gamma\|_{L^\infty})\,\|\gamma_{n-1}-\gamma\|_{L^1} \le 2\|\gamma-\gamma_{n-1}\|_{L^1} \to 0,
\]
and $\rho_{n-1}\to\rho_{\hc}$ (Lemma~\ref{lem:convergence}). Testing against $\varphi\in C_c^\infty(\{|x|>1\})$ and using dominated convergence on the left-hand side, we may pass to the limit:
\[
(-\Delta+4e)\gamma = 2e\rho_{\hc}(\gamma\ast\gamma) \quad\text{on } \{|x|>1\}, \qquad \gamma=1 \text{ on } \partial B_1.
\]
Since the right-hand side is bounded ($\gamma\ast\gamma\le\|\gamma\|_{L^1}$), elliptic regularity makes $\gamma$ continuous on $\{|x|>1\}$, and the equation holds classically. 
\newline

\textit{\textbf{Identification with $\tau$.}} Via the uniqueness property of $\tau$ from Lemma~\ref{lem:uniqueness}, it results immediately that $\gamma=\tau$.
%\textit{\textbf{Identification with $\tau$.}} We claim $\gamma\ge\tau_n$ for all $n$; for $n=0$ it is trivial. Assume $\gamma\ge\tau_{n-1}$. On $\{|x|>1\}$, subtracting \eqref{eq:iterate} from the limit equation,
%\[
%(-\Delta+4e)(\gamma-\tau_n) = 2e[\rho_{\hc}\,\gamma\ast\gamma - \rho_{n-1}\tau_{n-1}\ast\tau_{n-1}] \ge 0
%\]
%by the induction hypothesis and $\rho_{\hc}\ge\rho_{n-1}$ (Lemmas~\ref{lem:properties} and~\ref{lem:convergence}), while $\gamma-\tau_n=0$ on $\partial B_1$ (both equal $1$, and $\gamma$ is continuous up to the boundary by the first step) and $\liminf_{|x|\to\infty}(\gamma-\tau_n)\ge0$ since $\gamma\ge0$ and $\tau_n\to0$. Lemma~\ref{lem:maxp} with $D=\{|x|>1\}$, $c=4e$, $W=0$ gives $\gamma\ge\tau_n$; letting $n\to\infty$ yields $\gamma\ge\tau$, hence $\gamma=\tau$.
\end{proof}

\subsection{The $\omega_n$ iteration}
We recall the definition of the $\omega_n$ iteration
\begin{equation}\label{eq:omega2}
\omega_0=0
,\quad
(-\Delta + 4e + v_n)\omega_n = v_n + 2e\rho_{\hc}\omega_{n-1} * \omega_{n-1}.
\end{equation}
In this section we prove that $\omega_n$ converges to $\tau$ in $L^1$.

\begin{lemma}\label{lem:omega}
For all $n \ge 0$, $0 \leq \omega_n \leq 1$, $\omega_n \leq \tau$, $\omega_n \geq \gamma_n$, and $\omega_n \le \omega_{n+1}$.
\end{lemma}
\begin{proof}
Take $\omega_0=0$. By induction; the case $n=0$ is trivial. Assume $0\le\omega_{n-1}\le1$, $\omega_{n-1}\le\tau$, and $\omega_{n-1}\ge\gamma_{n-1}$.
\newline
\textit{\textbf{Non-negativity.}} The right-hand side of \eqref{eq:typeII} is non-negative and $(-\Delta+4e+v_n)^{-1}$ is positivity preserving, so $\omega_n\ge0$.
\newline
\textit{\textbf{Upper bound $\omega_n\le1$.}} On $\R^3$,
\[
(-\Delta+4e+v_n)(1-\omega_n) = 4e - 2e\rho_{\hc}\,\omega_{n-1}*\omega_{n-1} \ge 4e-2e = 2e > 0,
\]
using $\omega_{n-1}\le1$ and $\rho_{\hc}\int\omega_{n-1}\le\rho_{\hc}\int\tau=1$ (which follows from the inductive hypothesis and Lemma~\ref{lem:convergence}); moreover $1-\omega_n\to1$ at infinity. Lemma~\ref{lem:maxp} with $D=\R^3$, $c=4e$, $W=v_n$ gives $\omega_n\le1$.
\newline
\textit{\textbf{Monotonicity.}} Assuming $\omega_{n-1}\le\omega_n$, on $\R^3$,
\begin{align*}
&(-\Delta+4e+v_{n+1})(\omega_{n+1}-\omega_n) \\
&\quad = (v_{n+1}-v_n)(1-\omega_n) + 2e\rho_{\hc}[\omega_n*\omega_n - \omega_{n-1}*\omega_{n-1}] \ge 0,
\end{align*}
since $v_{n+1}\ge v_n$, $\omega_n\le1$, and $\omega_n\ge\omega_{n-1}\ge0$; moreover $\omega_{n+1}-\omega_n\to0$ at infinity. Lemma~\ref{lem:maxp} with $D=\R^3$, $c=4e$, $W=v_{n+1}$ gives $\omega_{n+1}\ge\omega_n$.
\newline
\textit{\textbf{Comparison $\omega_n\le\tau$.}} For $|x|\le 1$ this reads $\omega_n\le1=\tau$. On $\{|x|>1\}$, \eqref{eq:hsimpleq} gives $(-\Delta+4e)\tau = 2e\rho_{\hc}\,\tau*\tau$, so
\[
(-\Delta+4e)(\tau-\omega_n) = 2e\rho_{\hc}[\tau*\tau-\omega_{n-1}*\omega_{n-1}] \ge 0
\]
by the induction hypothesis, while $\tau-\omega_n = 1-\omega_n(1)\ge0$ on $\partial B_1$ and $\liminf_{|x|\to\infty}(\tau-\omega_n)\ge0$ since $\tau\ge0$ and $\omega_n\to0$. Lemma~\ref{lem:maxp} with $D=\{|x|>1\}$, $c=4e$, $W=0$ gives $\omega_n\le\tau$.
\newline
\textit{\textbf{Comparison $\omega_n\ge\gamma_n$.}} On $\R^3$, subtracting \eqref{eq:typeI} from \eqref{eq:typeII},
\[
(-\Delta+4e+v_n)(\omega_n-\gamma_n) = 2e[\rho_{\hc}\,\omega_{n-1}*\omega_{n-1} - \rho_{n-1}\gamma_{n-1}*\gamma_{n-1}] \ge 0,
\]
since $\rho_{\hc}\ge\rho_{n-1}$ (Lemmas~\ref{lem:properties} and~\ref{lem:convergence}) and $\omega_{n-1}\ge\gamma_{n-1}\ge0$ by the induction hypothesis and Lemma~\ref{lem:gamma-bounds}; moreover $\omega_n-\gamma_n\to0$ at infinity. Lemma~\ref{lem:maxp} with $D=\R^3$, $c=4e$, $W=v_n$ gives $\omega_n\ge\gamma_n$.
\end{proof}

\begin{corollary}\label{cor:omega-u}
$\|\omega_n - \tau\|_{L^1} \to 0$.
\end{corollary}

\begin{proof}
From $\gamma_n \le \omega_n \le \tau$ and Lemma~\ref{lem:gamma-limit} we have $\int(\tau-\omega_n) \le \int(\tau-\gamma_n) \to 0$, which implies $L^1$ convergence.
\end{proof}

%=====================================================================
\section{The Soft-Core Simple Equation}
%=====================================================================

%For $m > 0$, let $v_m = m\mathbbm{1}_{\Bone}$. The soft-core Simple Equation is
%\begin{equation}\label{eq:softcore}
%	(-\Delta + 4e + v_m)\sigma_{m,n} = v_m + 2e\bar{\rho}_{m,n-1}(\sigma_{m,n-1} * \sigma_{m,n-1}), \qquad
%	\bar{\rho}_{m,n} = \frac{2e}{m\int_{\Bone}(1 - \sigma_{m,n})}.
%\end{equation}
%The soft-core iterates $(\sigma_{m,n}, \bar{\rho}_{m,n})$ are defined by $\sigma_{m,0} = 0$, $\bar{\rho}_{m,0} = 0$, and
%\begin{equation}\label{eq:softcore-iterate}
%(-\Delta + 4e + v_m)\sigma_{m,n} = v_m + 2e\bar{\rho}_{m,n-1}\sigma_{m,n-1} * \sigma_{m,n-1},
%\end{equation}
%with $\bar{\rho}_{m,n} = 2e/(m\int_{\Bone}(1 - \sigma_{m,n}))$.

We recall \cite[Lemma 3.1-3.2, Equation (3.9)]{CJL20}:

\begin{lemma}\label{lem:softcore-props}
For all $n \geq 1$, $\sigma_{m,n} \geq 0$, $\sigma_{m,n} \leq 1$, $\sigma_{m,n} \leq \sigma_{m,n+1}$, $\bar{\rho}_{m,n} \leq \bar{\rho}_{m,n+1}$, and
\begin{equation}\label{eq:softcore-integral}
2\int\sigma_{m,n} = 1/\bar{\rho}_{m,n} + \bar{\rho}_{m,n-1}\|\sigma_{m,n-1}\|_{L^1}^2, \qquad \int\sigma_{m,n} \leq 1/\bar{\rho}_{m,n}.
\end{equation}
\end{lemma}

The proof uses the same tools as the proofs in section \ref{sec:intermediate}.
%\begin{remark}
%This follows from the analysis in \cite{CJL20}; we sketch the arguments for completeness. Non-negativity is due to positivity of the resolvent. Monotonicity in $n$ follows from the difference equation and the induction hypothesis. The upper bound $\sigma_{m,n}\le 1$ is obtained by a maximum principle argument similar to Lemma~\ref{lem:gamma-bounds}. The integral identity is obtained by integrating \eqref{eq:softcore-iterate} and using the definition of $\bar{\rho}_{m,n}$.
%\end{remark}

\subsection{Per-$n$ convergence to hard-core iterates}

\begin{lemma}\label{lem:per-n}
For each fixed $n \ge 1$, $\|\sigma_{m,n} - \tau_n\|_{L^1} \to 0$ and $\bar{\rho}_{m,n} \to \rho_n$ as $m \to \infty$.
\end{lemma}

\begin{proof}
By induction on $n$. For $n=1$, solving \eqref{eq:softcore} with $\sigma_{m,0}=0$ and the boundary conditions $\partial_r\sigma_{m,1}(r)|_{r=0}=0$ (which follows from the condition that $\sigma_{m,1}$ is differentiable and radially symmetric), $\lim_{r\to\infty}\sigma_{m,1}(r)=0$, along with the condition that $\sigma_{m,1}(r)$ and $\partial_r \sigma_{m,1}(r)$ are continuous at $r=1$, we find
\[
\sigma_{m,1}(r) = \begin{cases}
\dfrac{m}{\kappa^2}\Bigl(1 - \dfrac{(1+2\sqrt{e})}{\kappa\cosh\kappa+2\sqrt{e}\sinh\kappa}\,\dfrac{\sinh(\kappa r)}{r}\Bigr), & r<1,\\[8pt]
\sigma_{m,1}(1)\,\dfrac{e^{-2\sqrt{e}(r-1)}}{r}, & r>1,
\end{cases}
\qquad \kappa = \sqrt{4e+m},
\]
with $\sigma_{m,1}(1) = \frac{m}{\kappa^2}\,\frac{\kappa\cosh\kappa-\sinh\kappa}{\kappa\cosh\kappa+2\sqrt{e}\sinh\kappa} = 1-O(1/\sqrt{m})$.
In addition, it follows from integrating the first of \eqref{eq:softcore} that
\begin{equation}
  m\int_{\Bone}(1-\sigma_{m,1}(|x|))\ dx = 4e\int_{\mathbb R^3}\sigma_{m,1}(|x|)\ dx
  \label{integrated}
\end{equation}
so in particular $\int_{\Bone}|\sigma_{m,1}-\tau_1|\equiv \int_{\Bone}(1-\sigma_{m,1})=O(1/m)$.
Combining this with $\sigma_{m,1}(1)=1+O(1/\sqrt m)$ and \eqref{eq:u1}, we find that $\|\sigma_{m,1}-\tau_1\|_{L^1} = O(1/\sqrt{m})$.
Finally, by \eqref{integrated}, $m\int_{\Bone}(1-\sigma_{m,1}) = 4e\int\sigma_{m,1} \to 4e\int\tau_1 = 2e/\rho_1$ (by Eq.~\eqref{eq:u1}), and therefore by the second of \eqref{eq:softcore}, $\bar{\rho}_{m,1}\to\rho_1$.
\newline 
\textit{\textbf{Estimates inside the core.}} We mimic the \textit{Boundary value} step in the proof of Lemma~\ref{lem:gamma-limit}, with the roles of $n$ there played by $m$. Set $\iota_{m,n}= 1-\sigma_{m,n} \ge 0$ on $\Bone$; it satisfies, for $x\in\Bone$,
\[
	(-\Delta+\kappa^2)\iota_{m,n} = 4e - 2e\bar{\rho}_{m,n-1}(\sigma_{m,n-1}\ast\sigma_{m,n-1}) \equiv f,
\qquad \kappa=\sqrt{4e+m},
\]
with $0\le f\le 4e$: indeed $2e\bar{\rho}_{m,n-1}(\sigma_{m,n-1}\ast\sigma_{m,n-1}) \le 2e\bar{\rho}_{m,n-1}\|\sigma_{m,n-1}\|_{L^\infty}\|\sigma_{m,n-1}\|_{L^1} \le 2e$ by Lemma~\ref{lem:softcore-props}. Decompose $\iota_{m,n}  = \tilde\iota_{m,n}  + h_{m,n}$, where $\tilde\iota_{m,n} |_{\partial B_1}=0$ and $h_{m,n}$ is the homogeneous solution: $(-\Delta+\kappa^2)h_{m,n}=0$ with $h_{m,n}(1)=\iota_{m,n}(1)$, i.e.\ $h_{m,n}(r) = \iota_{m,n}(1)\frac{\sinh(\kappa r)}{r\sinh\kappa}$. Since $(-\Delta+\kappa^2)h_{m,n}=0$, the function $\tilde\iota_{m,n} $ satisfies the same equation as $\iota_{m,n}$, with $\tilde\iota_{m,n} =0$ on $\partial B_1$; Corollary~\ref{cor:maxp} with $D=\Bone$, $c=\kappa^2$, $W=0$, $F=4e$ gives $0\le\tilde\iota_{m,n} \le4e/\kappa^2$. Furthermore,
\begin{equation}
	\int_{\Bone}h_{m,n} = \frac{4\pi\iota_{m,n} (1)}{\kappa^2\sinh\kappa}(\kappa\cosh\kappa-\sinh\kappa) \sim \frac{4\pi\iota_{m,n} (1)}{\kappa}\quad(\kappa\to\infty).
	\label{inthmn}
\end{equation}
Integrating \eqref{eq:softcore} over $\R^3$ and using Lemma~\ref{lem:softcore-props},
\[
m\int_{\Bone}\iota_{m,n}  = 4e\int\sigma_{m,n} - 2e\bar{\rho}_{m,n-1}\|\sigma_{m,n-1}\|_{L^1}^2 \le \frac{4e}{\bar{\rho}_{m,1}},
\]
and since $\bar{\rho}_{m,1}\to\rho_1>0$, the right-hand side is bounded uniformly in $m$, so $\int_{\Bone}\iota_{m,n}  = O(1/m)$. Since $\iota_{m,n} = \tilde\iota_{m,n}  + h_{m,n} \ge h_{m,n}$, we thus have $\int_{\Bone}h_{m,n}=O(1/m)$, which, by \eqref{inthmn}, implies $\iota_{m,n} (1) = O(\kappa/m) = O(1/\sqrt{m})$. Hence $1-\sigma_{m,n}(1) = O(m^{-1/2})$ and $\int_{\Bone}(1-\sigma_{m,n}) = O(m^{-1})$.
\newline
%\textit{\textbf{Convergence outside $\Bone$.}} Let $w = \tau_n - \sigma_{m,n}$ on $\{|x|>1\}$. Subtracting \eqref{eq:softcore} from \eqref{eq:iterate},
%\[
%(-\Delta+4e)w = 2e\bigl(\rho_{n-1}\tau_{n-1}\ast\tau_{n-1} - \bar{\rho}_{m,n-1}\sigma_{m,n-1}\ast\sigma_{m,n-1}\bigr) \equiv \beta_n.
%\]
%By the Green's function representation (cf.\ \eqref{eq:un-exterior}, with $\mu=2\sqrt{e}$), {\color{red}Define $w$ on $\mathbb R^3$ as $w=\tau_n-\sigma_{m,n}$ for $|x|>1$ and $(-\Delta+4e)w=\beta_n$ for $|x|\le 1$. In this case, for $|x|>1$,
%\[
%w(r) = b Y_{4e} + Y_{4e}\ast \beta_n.
%\]
%}
%\[
%w(r) = b\,\frac{e^{-2\sqrt{e}(r-1)}}{r} + \frac{1}{4\sqrt{e}\,r}\int_1^\infty s\,\beta_n(s)\,e^{-2\sqrt{e}|s-r|}\,ds.
%\]
%\[
%  b=w(1)- \frac{1}{4\sqrt{e}}\int_1^\infty s\,\beta_n(s)\,e^{-2\sqrt{e}(s-1)}\,ds
%\]
%Taking $L^1$ norms, using $\int_1^\infty 4\pi r\,e^{-2\sqrt{e}(r-1)}\,dr = \frac{4\pi(1+2\sqrt{e})}{4e}$ for the first term and Young's inequality (the Yukawa kernel has $\|Y_{4e}\|_{L^1}=1/(4e)$) for the second,
%\[
%\|w\|_{L^1(B_1^c)} \le |w(1)|\,\frac{4\pi(1+2\sqrt{e})}{4e} + \frac{\|\beta_n\|_{L^1}}{4e}.
%\]

\textit{\textbf{Convergence outside $\Bone$.}} Let $w = \tau_n - \sigma_{m,n}$ on $\{|x|>1\}$; $w$ is radial and bounded, with $w(1) = 1-\sigma_{m,n}(1)$. Subtracting \eqref{eq:softcore} from \eqref{eq:iterate},
\[
(-\Delta+4e)w = 2e\bigl(\rho_{n-1}\tau_{n-1}\ast\tau_{n-1} - \bar{\rho}_{m,n-1}\sigma_{m,n-1}\ast\sigma_{m,n-1}\bigr) \equiv \beta_n \quad\text{on }\{|x|>1\}.
\]
Since $w$ solves this equation only on the exterior, it cannot be represented directly by the Green's function $Y_{4e}(x) = e^{-2\sqrt{e}|x|}/(4\pi|x|)$ of $-\Delta+4e$ on $\R^3$. Instead, we extend the source term $\beta_n$ to $\mathbb R^n$ by setting it to zero inside $\Bone$: let $W := Y_{4e}\ast(\beta_n\mathbbm{1}_{B_1^c})$ be the bounded radial solution of $(-\Delta+4e)W = \beta_n\mathbbm{1}_{B_1^c}$ on $\R^3$.
On $\{|x|>1\}$ it is given by directly evaluating the convolution of $\beta_n$ with $Y_{4e}$ and using bipolar coordinates:
\[
W(r) = \frac{1}{4\sqrt{e}\,r}\int_1^\infty s\,\beta_n(s)\bigl(e^{-2\sqrt{e}|r-s|} - e^{-2\sqrt{e}(r+s)}\bigr)\,ds
\]
so it is $C^2$ there and solves the equation classically.
Both $w$ and $W$ solve $(-\Delta+4e)u=\beta_n$ on the exterior, and we set $W(1)=w(1)$ so that $w(r)=W(r)$ for all $r \ge 1$, and so, for $r \ge 1$,
\begin{equation}
  w(r)=Y_{4e}\ast(\beta_n \mathbbm 1_{\Bone^c})
  .
\end{equation}
We then estimate $\|w\|_{L^1(B_1^c)}\le\|w\|_{L^1}$, and, by Young's inequality and $\|Y_{4e}\|_{L^1} = 1/(4e)$,
\[
\|w\|_{L^1(B_1^c)} \le \frac{\|\beta_n\|_{L^1}}{4e}
\]
By the induction hypothesis and Young's inequality,
\[
\|\beta_n\|_{L^1} \le 2e|\rho_{n-1}-\bar{\rho}_{m,n-1}|\,\|\tau_{n-1}\|_{L^1}^2 + 2e\bar{\rho}_{m,n-1}\bigl(\|\tau_{n-1}\|_{L^1}+\|\sigma_{m,n-1}\|_{L^1}\bigr)\|\tau_{n-1}-\sigma_{m,n-1}\|_{L^1} \to 0
.
\]
Combining with $\int_{\Bone}|\tau_n-\sigma_{m,n}| \equiv \int_{\Bone}(1-\sigma_{m,n}) = O(1/m)$, we conclude $\|\sigma_{m,n}-\tau_n\|_{L^1} \to 0$.
\newline

%The first term is $O(1/\sqrt{m})$ by the core estimate ($w(1)=1-\sigma_{m,n}(1)\ge0$). For the second,
%\[
%\|\beta_n\|_{L^1} \le 2e|\rho_{n-1}-\bar{\rho}_{m,n-1}|\,\|\tau_{n-1}\|_{L^1}^2 + 2e\bar{\rho}_{m,n-1}\bigl(\|\tau_{n-1}\|_{L^1}+\|\sigma_{m,n-1}\|_{L^1}\bigr)\|\tau_{n-1}-\sigma_{m,n-1}\|_{L^1} \to 0
%\]
%by the induction hypothesis. Combining with $\int_{\Bone}|\tau_n-\sigma_{m,n}| = \int_{\Bone}(1-\sigma_{m,n}) = O(1/m)$, we conclude $\|\sigma_{m,n}-\tau_n\|_{L^1} \to 0$.

\textit{\textbf{Density convergence.}} By the integral relation \eqref{eq:softcore-integral} and the $L^1$ convergence just proved,
\[
\frac{1}{\bar{\rho}_{m,n}} = 2\int\sigma_{m,n} - \bar{\rho}_{m,n-1}\Bigl(\int\sigma_{m,n-1}\Bigr)^2 \;\longrightarrow\; 2\int\tau_n - \rho_{n-1}\Bigl(\int\tau_{n-1}\Bigr)^2 = \frac{1}{\rho_n},
\]
the last equality being \eqref{eq:balance}; hence $\bar{\rho}_{m,n}\to\rho_n$.
\end{proof}

\section{Convergence of the Converged Soft-Core Solution}
%=====================================================================

Let $\sigma_m = \lim_{n\to\infty}\sigma_{m,n}$ and $\bar{\rho}_m = \lim_{n\to\infty}\bar{\rho}_{m,n}$. Then, by \cite[Lemma 3.3]{CJL20}, $(\sigma_m,\bar{\rho}_m)$ solves\
\begin{equation}
  (-\Delta+4e+v_m)\sigma_m=v_m+2e\bar\rho_m \sigma_m\ast \sigma_m
\end{equation}
with
\begin{equation}
  \bar\rho_m=\frac{2e}{m\int_{\Bone}(1-\sigma_m)}
  .
  \label{eq:softcore_lim}
\end{equation}
In particular, we have $\bar{\rho}_m\int\sigma_m = 1$.

\begin{lemma}\label{lem:rhom-decreasing}
The map $m\mapsto\bar{\rho}_m$ is differentiable and strictly decreasing; consequently $\bar{\rho}_m \ge \rho_{\hc}$ for all $m$.
\end{lemma}

\begin{proof}
\textit{\textbf{The derivative formula.}}
For $f\in L^1$, let $C_f$ denote the convolution operator $C_f\,g := f*g$, and define the linearised resolvent
\[
K_m = \bigl(-\Delta+4e+v_m - 4e\bar{\rho}_m C_{\sigma_m}\bigr)^{-1}.
\]
%$K_m$ is self-adjoint (as $v_m$ is real and $\sigma_m$ is radial, so $C_{\sigma_m}$ is self-adjoint), positivity preserving by the Trotter product formula, and $\Phi_m := K_m[v_m]$ satisfies $0\le\Phi_m\le1$, with $\Phi_m<1$ on a set of positive measure \cite[Lemma~7.5]{CJL21}. 
The operator $K_m^{-1}$ is symmetric and non-negative on $L^2$ (as $v_m$ is real and $\sigma_m$ is radial, so $C_{\sigma_m}$ is self-adjoint with norm $\|\sigma_m\|_{L^1}=1/\bar{\rho}_m$); $K_m[v_m] := \int_0^\infty e^{-tK_m^{-1}}v_m\,dt$ is well defined and non-negative by the Trotter product formula, and $\Phi_m := K_m[v_m]$ satisfies $0\le\Phi_m\le1$, with $\Phi_m<1$ on a set of positive measure \cite[Lemma~7.5]{CJL21}.

Since $\bar{\rho}_m\int\sigma_m=1$,
\[
K_m^{-1}\,1 = 4e+v_m-4e\bar{\rho}_m\!\int\!\sigma_m = v_m = K_m^{-1}\,\Phi_m,
\]
so $1-\Phi_m$ solves the homogeneous equation
\begin{equation}\label{eq:hom-phi}
(-\Delta+4e+v_m)(1-\Phi_m) = 4e\bar{\rho}_m\,\sigma_m*(1-\Phi_m).
\end{equation}
Multiplying \eqref{eq:softcore_lim} by $1-\Phi_m$ and \eqref{eq:hom-phi} by $\sigma_m$, and integrating over $\R^3$ gives
\begin{align}
\int(1-\Phi_m)\,(-\Delta+4e+v_m)\sigma_m &= \int(1-\Phi_m)\,v_m + 2e\bar{\rho}_m\int(1-\Phi_m)\,(\sigma_m\ast\sigma_m), \label{eq:pairA}\\
\int\sigma_m\,(-\Delta+4e+v_m)(1-\Phi_m) &= 4e\bar{\rho}_m\int\sigma_m\,\bigl(\sigma_m\ast(1-\Phi_m)\bigr). \label{eq:pairB}
\end{align}
All integrals converge: the right-hand side of \eqref{eq:softcore_lim} lies in $L^1$ and $0\le1-\Phi_m\le1$, while $\sigma_m\in L^1$ and $\sigma_m\ast(1-\Phi_m)$ is bounded. We prove that the two left-hand sides coincide.
The terms with $4e+v_m$ agree trivially.
For the Laplacian terms we use integration by parts.
To show that the boundary terms vanish, we use the fact that $\sigma_m(x)=O(|x|^{-4})$ \cite{CJL20} along with
\begin{equation}
  \nabla\sigma_m = \nabla Y_{4e}\ast (v_m(1-\sigma_m)+2e\bar{\rho}_m\,\sigma_m\ast\sigma_m) = O(|x|^{-4})
\end{equation}
while $|\Phi_m|\le1$ and
\begin{equation}
  \nabla\Phi_m = \nabla Y_{4e}\ast (v_m(1-\Phi_m)+4e\bar{\rho}_m\,\sigma_m\ast\Phi_m)=O(1)
  .
\end{equation}
Thus,
\begin{equation}
  \int(1-\Phi_m)\,v_m + 2e\bar{\rho}_m\int(1-\Phi_m)\,(\sigma_m\ast\sigma_m) =
  4e\bar{\rho}_m\int\sigma_m\,\bigl(\sigma_m\ast(1-\Phi_m)\bigr).
\end{equation}
Furthermore, Fubini's theorem and the evenness of the radial function $\sigma_m$ give
\[
\int\sigma_m(x)\,\bigl(\sigma_m\ast(1-\Phi_m)\bigr)(x)\,dx
= \int(1-\Phi_m)\,(\sigma_m\ast\sigma_m)
\]
and so
\begin{equation}\label{eq:R}
\int(1-\Phi_m)\,v_m = 2e\bar{\rho}_m\int(1-\Phi_m)\,(\sigma_m\ast\sigma_m).
\end{equation}

Now differentiate in $m$. To justify this we argue with finite differences, following \cite[Section~3]{CJL21}. Fix $m'\neq m$ and write $\delta m := m'-m$, $\delta\sigma := \sigma_{m'}-\sigma_m$, $\delta\rho := \bar{\rho}_{m'}-\bar{\rho}_m$, $\tilde\sigma := \tfrac12(\sigma_m+\sigma_{m'})$, and
\[
\tilde K^{-1} := -\Delta+4e+v_{m'}-4e\bar{\rho}_{m'}C_{\tilde\sigma};
\]
by \cite[Section~3]{CJL21} there is a bounded radial function $\tilde\Phi$ with $\tilde K^{-1}\tilde\Phi = v_{m'}$ and $0\le\tilde\Phi\le1$, and, as for $\Phi_m$, $\nabla\tilde\Phi$ is bounded. Subtracting \eqref{eq:softcore_lim} at $m$ from \eqref{eq:softcore_lim} at $m'$, and using $\sigma_{m'}\ast\sigma_{m'}-\sigma_m\ast\sigma_m = 2\,\delta\sigma\ast\tilde\sigma$ and $v_{m'}-v_m = \delta m\,\mathbbm{1}_{\Bone}$, gives
\begin{equation}\label{eq:sigma-prime}
\tilde K^{-1}\delta\sigma = \delta m\,\mathbbm{1}_{\Bone}(1-\sigma_m) + 2e\,\delta\rho\,(\sigma_m*\sigma_m),
\end{equation}
while subtracting $2e/\bar{\rho}_m = m\int_{\Bone}(1-\sigma_m)$ from $2e/\bar{\rho}_{m'} = m'\int_{\Bone}(1-\sigma_{m'})$ gives
\begin{equation}
-\frac{2e\,\delta\rho}{\bar{\rho}_m\bar{\rho}_{m'}} = \delta m\int_{\Bone}(1-\sigma_m) - \int v_{m'}\,\delta\sigma.
\label{rho'intvs}
\end{equation}
Multiplying $\tilde K^{-1}\tilde\Phi = v_{m'}$ by $\delta\sigma$ and integrating, and moving the operator onto $\delta\sigma$ --- the Laplacian terms by integration by parts, whose boundary terms vanish as above since $\delta\sigma,\nabla\delta\sigma = O(|x|^{-4})$ while $\tilde\Phi,\nabla\tilde\Phi$ are bounded, and the convolution terms by Fubini's theorem and the evenness of $\tilde\sigma$ --- we find $\int v_{m'}\,\delta\sigma = \int\tilde\Phi\,\tilde K^{-1}\delta\sigma$, hence by \eqref{eq:sigma-prime}
\begin{equation}
\int v_{m'}\,\delta\sigma = \int \tilde\Phi\Bigl[\delta m\,\mathbbm{1}_{\Bone}(1-\sigma_m) + 2e\,\delta\rho\,(\sigma_m*\sigma_m)\Bigr].
\label{intvs}
\end{equation}
Using $\int(\sigma_m*\sigma_m) = \bigl(\int\sigma_m\bigr)^2 = 1/\bar{\rho}_m^2$ and $\delta\rho=\bar\rho_{m'}-\bar\rho_m$, we write
\begin{equation}
  -2e\,\delta\rho\int \sigma_m\ast \sigma_m + \frac{2e\,(\delta\rho)^2}{\bar{\rho}_m^2\bar{\rho}_{m'}}
  =
  -\frac{2e\,\delta\rho}{\bar\rho_m^2} + \frac{2e\,\delta\rho(\bar\rho_{m'}-\bar\rho_m)}{\bar{\rho}_m^2\bar{\rho}_{m'}}
  =
  -\frac{2e\,\delta\rho}{\bar{\rho}_m\bar{\rho}_{m'}}
  \label{rho'}
\end{equation}
so, inserting \eqref{intvs} and \eqref{rho'} into \eqref{rho'intvs}, we find
\[
0 = \delta m\int_{\Bone}(1-\tilde\Phi)(1-\sigma_m) + 2e\,\delta\rho\Bigl[\int(1-\tilde\Phi)(\sigma_m*\sigma_m) - \frac{\delta\rho}{\bar{\rho}_m^2\bar{\rho}_{m'}}\Bigr].
\]
Solving for the difference quotient, we arrive at the exact identity
\begin{equation}\label{eq:fd-quotient}
\frac{\delta\rho}{\delta m}\equiv\frac{\bar{\rho}_{m'}-\bar{\rho}_m}{m'-m} = -\,\frac{\int_{\Bone}(1-\tilde\Phi)(1-\sigma_m)}{2e\Bigl[\int(1-\tilde\Phi)\,(\sigma_m*\sigma_m) - \dfrac{\delta\rho}{\bar{\rho}_m^2\bar{\rho}_{m'}}\Bigr]},
.
\end{equation}
Now, following \cite[Section~3 and Lemma~7.7]{CJL21}, we can show that as $m'\to m$, $\delta\rho\to0$ and $\tilde\Phi\to\Phi_m$ almost everywhere; (the arguments in \cite[Section~3 and Lemma~7.7]{CJL21} apply almost exactly using the fact that  the dependence of $K^{-1}$ on $m$ is linear). Since $0\le\tilde\Phi\le1$, dominated convergence passes \eqref{eq:fd-quotient} to the limit, and eliminating the convolution integral via \eqref{eq:R}, which gives $2e\int(1-\Phi_m)(\sigma_m*\sigma_m) = \frac{m}{\bar{\rho}_m}\int_{\Bone}(1-\Phi_m)$, we arrive at
\begin{equation}\label{eq:rhoderiv}
\partial_m\bar{\rho}_m = -\frac{\bar{\rho}_m}{m}\,
\frac{\int_{\Bone}(1-\Phi_m)(1-\sigma_m)}{\int_{\Bone}(1-\Phi_m)}.
\end{equation}
We can now prove that $\bar\rho$ is strictly decreasing: $1-\Phi_m\ge0$ does not vanish identically, so by \eqref{eq:hom-phi} and the strict positivity of the kernel of $(-\Delta+4e+v_m)^{-1}$, in fact $1-\Phi_m>0$ everywhere; since moreover $\int_{\Bone}(1-\sigma_m) = 2e/(m\bar{\rho}_m)>0$, both the numerator and the denominator of \eqref{eq:rhoderiv} are strictly positive, whence $\partial_m\bar{\rho}_m<0$.
\newline
\textit{\textbf{Comparison with the hard-core density.}} By Lemma~\ref{lem:softcore-props} the iterates $\sigma_{m,n}$ increase in $n$, so $\bar{\rho}_{m,n} = 2e/(m\int_{\Bone}(1-\sigma_{m,n}))$ increases in $n$ towards $\bar{\rho}_m$; in particular $\bar{\rho}_m\ge\bar{\rho}_{m,n}$ for every $n$. By Lemma~\ref{lem:per-n}, $\bar{\rho}_{m,n}\to\rho_n$ as $m\to\infty$, hence $\liminf_{m\to\infty}\bar{\rho}_m\ge\rho_n$ for every $n$, and letting $n\to\infty$, $\liminf_{m\to\infty}\bar{\rho}_m\ge\rho_{\hc}$. Since $\bar{\rho}_m$ is decreasing in $m$, the limit $\lim_{m\to\infty}\bar{\rho}_m$ exists and is $\ge\rho_{\hc}$; monotonicity then gives $\bar{\rho}_m\ge\rho_{\hc}$ for every $m$.
\end{proof}

\begin{lemma}\label{lem:comparison}
For all $m \ge 0 $, $\sigma_m \ge \omega_m$. 
\end{lemma}

\begin{proof}
Induction on $m$. The case $m=0$ is trivial. Assuming $\sigma_{m-1} \ge \omega_{m-1}$ and recalling that $v_m = m\mathbf 1_{\Bone}$
\begin{align*}
&(-\Delta+4e+v_m)(\sigma_m-\omega_m) = 2e\bigl[\bar{\rho}_{m-1}\sigma_{m-1}\ast\sigma_{m-1} - \rho_{\hc}\omega_{m-1}\ast\omega_{m-1}\bigr] \ge 0,
\end{align*}
because $\bar{\rho}_{m-1}\ge\rho_{\hc}$ and $\sigma_{m-1}\ge\omega_{m-1}$. Positivity-preservation of $(-\Delta+4e+v_m)^{-1}$ yields the lemma.
\end{proof}

\begin{theorem}\label{thm:convergencesigma}
$\|\sigma_m - \tau\|_{L^1} \to 0$ and $\bar{\rho}_m \to \rho_{\hc}$ as $m \to \infty$.
\end{theorem}

\begin{proof}
From Lemma~\ref{lem:comparison} we have $\sigma_m \ge \omega_m$. Thus
\[
\frac{1}{\bar{\rho}_m} = \int\sigma_m \ge \int\omega_m \longrightarrow \frac{1}{\rho_{\hc}},
\]
by Corollary~\ref{cor:omega-u}. Hence $\limsup_{m\to\infty}\bar{\rho}_m \le \rho_{\hc}$, which together with Lemma~\ref{lem:rhom-decreasing} gives $\bar{\rho}_m\to\rho_{\hc}$.
For the $L^1$ convergence,
\[
\|\sigma_m - \tau\|_{L^1} \le \|\sigma_m - \omega_m\|_{L^1} + \|\omega_m - \tau\|_{L^1}.
\]
The second term tends to zero by Corollary~\ref{cor:omega-u}. Since $\sigma_m \ge \omega_m$,
\[
\|\sigma_m - \omega_m\|_{L^1} = \int\sigma_m - \int\omega_m \to 0.
\]
Hence $\|\sigma_m - \tau\|_{L^1} \to 0$.
\end{proof}

%=====================================================================
\section{Rate of Convergence}
%=====================================================================

\begin{theorem}[Proof of Theorem~\ref{thm:OE}(i)]\label{thm:rate}
$|\bar{\rho}_m - \rho_{\hc}| = O(m^{-1/2})$.
\end{theorem}

\begin{proof}
From \eqref{eq:rhoderiv}, since the ratio appearing there is an average of $1-\sigma_m$ against the positive weight $(1-\Phi_m)\mathbbm{1}_{\Bone}$,
\[
|\partial_m\bar{\rho}_m| \le \frac{\bar{\rho}_m}{m}\,\sup_{\Bone}(1-\sigma_m).
\]
We estimate the supremum as in the core estimate of Lemma~\ref{lem:per-n}, now at the level of the converged solution. Set $\iota_{m} = 1-\sigma_m \ge 0$; on $\Bone$, \eqref{eq:softcore} at convergence reads
\[
(-\Delta+\kappa^2)\iota_m = 4e - 2e\bar{\rho}_m(\sigma_m*\sigma_m) \equiv f,
\qquad \kappa = \sqrt{4e+m},
\]
with $0\le f\le 4e$, since $\sigma_m\le1$ and $\bar{\rho}_m\int\sigma_m=1$ give $\bar{\rho}_m(\sigma_m*\sigma_m)\le1$. Decompose $\iota_m = \tilde\iota_m + h_m$, where $\tilde\iota_m|_{\partial B_1}=0$ and $h_m(r) = \iota_m(1)\frac{\sinh(\kappa r)}{r\sinh\kappa}$ is the homogeneous solution with $h_m(1)=\iota_m(1)$. Corollary~\ref{cor:maxp} with $D=\Bone$, $c=\kappa^2$, $W=0$, $F=4e$ gives $0\le\tilde\iota_m\le 4e/\kappa^2$, while $h_m\le\iota_m(1)$ on $\Bone$ because $r\mapsto\sinh(\kappa r)/r$ is increasing; hence
\[
\sup_{\Bone}(1-\sigma_m) \le (1-\sigma_m(1)) + \frac{4e}{\kappa^2}.
\]
For the boundary value, the definition of $\bar{\rho}_m$ gives exactly $\int_{\Bone}\iota_m = 2e/(m\bar{\rho}_m) \le 2e/(m\rho_{\hc})$, using $\bar{\rho}_m\ge\rho_{\hc}$ from Lemma~\ref{lem:rhom-decreasing}, while $\int_{\Bone}h_m \sim 4\pi\iota_m(1)/\kappa$ as $\kappa\to\infty$, as in Lemma~\ref{lem:per-n}. Since $\iota_m = \tilde\iota_m + h_m \ge h_m \ge 0$, this forces $\iota_m(1) = O(\kappa/m) = O(m^{-1/2})$, and altogether
\[
\sup_{\Bone}(1-\sigma_m) = O(m^{-1/2}).
\]
%Since $\bar{\rho}_m$ is decreasing, $\bar{\rho}_m\le\bar{\rho}_1$ for $m\ge1$, so $|\partial_m\bar{\rho}_m| = O(m^{-3/2})$. Finally, $\bar{\rho}_{m'}\to\rho_{\hc}$ as $m'\to\infty$ by Theorem~\ref{thm:convergencesigma}, so
%\[
%|\bar{\rho}_m - \rho_{\hc}| = \Bigl|\int_m^\infty \partial_{m'}\bar{\rho}_{m'}\,dm'\Bigr| \le C\int_m^\infty (m')^{-3/2}\,dm' = O(m^{-1/2}).
%\]
Since $\bar{\rho}_m$ is decreasing, $\bar{\rho}_m\le\bar{\rho}_1$ for $m\ge1$, so $|\partial_m\bar{\rho}_m| \le C m^{-3/2}$. For any $M>m$, the fundamental theorem of calculus gives $\bar{\rho}_M-\bar{\rho}_m = \int_m^M\partial_{m'}\bar{\rho}_{m'}\,dm'$; letting $M\to\infty$, the left side tends to $\rho_{\hc}-\bar{\rho}_m$ by Theorem~\ref{thm:convergencesigma}, while the right side converges absolutely by the bound above. Hence
\[
|\bar{\rho}_m - \rho_{\hc}| \le \int_m^\infty C\,(m')^{-3/2}\,dm' = \frac{2C}{\sqrt{m}}.
\]

\end{proof}

%\begin{proof}
%From the derivative formula in Lemma~\ref{lem:rhom-decreasing}, {\color{red}reference to the equation}
%\[
%|\partial_m\bar{\rho}_m| \le \frac{\bar{\rho}_m}{m}\,\sup_{\Bone}(1-\sigma_m).
%\]
%On $\Bone$, $1-\sigma_m$ satisfies the same equation as $\phi$ in Lemma~\ref{lem:per-n} with $\kappa=\sqrt{4e+m}$ {\color{red}reference, and write down the equation}. Decomposing into particular and homogeneous parts gives {\color{red}more details}
%\[
%\sup_{\Bone}(1-\sigma_m) \le (1-\sigma_m(1)) + \frac{4e}{\kappa^2} = O(m^{-1/2}).
%\]
%Hence $|\partial_m\bar{\rho}_m| = O(m^{-3/2})$. Integrating from $m$ to $\infty$ and using $\bar{\rho}_m\to\rho_{\hc}$ yields $|\bar{\rho}_m - \rho_{\hc}| = O(m^{-1/2})$.
%\end{proof}

%=====================================================================
\section{The LHY Formula}\label{sec:LHY}
The scattering solution $\phi_m$ of $v_m$ is defined by
\[
(-\Delta+v_m)\phi_m = v_m, \qquad \phi_m(x)\to0 \text{ as } |x|\to\infty;
\]
equivalently, $\psi_m := 1-\phi_m$ is the zero-energy scattering solution, $-\Delta\psi_m + v_m\psi_m = 0$ with $\psi_m\to1$ at infinity. For $v_m = m\mathbbm{1}_{\Bone}$ it is explicitly done by using the boundary conditions that require $\partial \psi_m (0)=0$ (due to differentiability and radial symmetry), and that $\psi_m$ and $\partial_m \psi$ be continuous at the boundary of the unit ball. This gives $\psi_m(r) = \frac{\sinh(\sqrt{m}\,r)}{r\sqrt{m}\cosh\sqrt{m}}$ for $r\le1$ and $\psi_m(r) = 1 - a_0(m)/r$ for $r\ge1$, where matching at $r=1$ gives the scattering length
\begin{equation}
	a_0(m) = 1-\frac{\tanh\sqrt{m}}{\sqrt{m}} \label{eq:tanh}
\end{equation}
which increases to $1$ as $m\to\infty$; moreover $0\le\phi_m\le1$. 
%We stress that $\phi_m$ is unrelated to the boundary-layer functions $\phi$, $\tilde\phi$, $\phi_n$ appearing in the proofs of Lemmas~\ref{lem:gamma-limit} and~\ref{lem:per-n}; we retain the symbol as it is standard for the scattering solution.

We adopt the Fourier transform convention 
\[
\hat f(k) = \int e^{ik\cdot x}f(x)\,dx, \qquad f(x) = \frac{1}{(2\pi)^3}\int e^{-ik\cdot x}\hat f(k)\,dk.
\]
For each $m>0$, define
\begin{equation}\label{eq:Sm}
S_m(k) = \frac{\bar{\rho}_m}{2e}\int e^{ik\cdot x}(1-\sigma_m(x))\,v_m(x)\,dx,
\end{equation}
and let $w_m$ be the decaying solution of $-\Delta w_m = (1-\sigma_m)v_m$, that is, $\hat w_m(k) = 2eS_m(k)/(\bar{\rho}_m k^2)$. Fourier transforming \eqref{eq:softcore} at convergence gives the quadratic equation $2e\bar{\rho}_m\hat\sigma_m^2 - (k^2+4e)\hat\sigma_m + \frac{2e}{\bar{\rho}_m}S_m = 0$; selecting the root that vanishes as $|k|\to\infty$ yields \cite[Eq.~(1.43)]{CJL21}
\begin{equation}\label{eq:eta-hat}
\hat\sigma_m(k) = \frac{1}{\bar{\rho}_m}\left(\frac{k^2}{4e}+1-\sqrt{\Bigl(\frac{k^2}{4e}+1\Bigr)^2 - S_m(k)}\right).
\end{equation}
(As a consistency check, note that $\hat\sigma_m(0)=1/\bar{\rho}_m$, which is the normalization $\bar{\rho}_m\int\sigma_m=1$.)

\begin{lemma}\label{lem:Sm-bound}
$S_m(k)$ is real, $|S_m(k)|\le1$ for all $k$, $m$, and $S_m(0)=1$.
\end{lemma}
\begin{proof}
$(1-\sigma_m)v_m$ is non-negative and radial, so $S_m$ is real with $|S_m(k)|\le S_m(0) = \frac{\bar{\rho}_m}{2e}\int(1-\sigma_m)v_m = 1$, the last equality follows from the definition of $\bar{\rho}_m$ (cf.\ \cite[Sec.~5]{CJL21}).
\end{proof}

\begin{lemma}\label{lem:eta-w-bound}
Let $g(\tilde k,s) := \tilde k^2+1-\sqrt{(\tilde k^2+1)^2-s}-\frac{s}{2\tilde k^2}$ and $h(\tilde k) := g(\tilde k,-1)-g(\tilde k,1)$. Then $h\in L^1(\R^3)$, and for all $x\in\R^3$, $m>0$,
\[
	|\sigma_m(x)-w_m(x)| \le \frac{C_0\,e^{3/2}}{\bar{\rho}_m}, \qquad C_0 := \frac{1}{\pi^3}\int h(\tilde k)\,d^3\tilde k.
\]
\end{lemma}
\begin{proof}
	Set $\tilde k = k/(2\sqrt e)$. By \eqref{eq:eta-hat}, the definition of $\hat w_m$,
\begin{equation}\label{eq:sigma-w-rep}
\sigma_m(x)-w_m(x) = \frac{e^{3/2}}{\bar{\rho}_m\pi^3}\int e^{-2i\sqrt e\,\tilde k\cdot x}\,g\bigl(\tilde k, S_m(2\sqrt e\,\tilde k)\bigr)\,d^3\tilde k
\end{equation}
(cf.\ \cite[Eq.~(5.13)]{CJL21} for the analogous representation). The map $s\mapsto g(\tilde k,s)$ is decreasing on $[-1,1]$, since $\partial_s g = \frac{1}{2\sqrt{(\tilde k^2+1)^2-s}}-\frac{1}{2\tilde k^2}\le0$, with $g(\tilde k,0)=0$; hence, by Lemma~\ref{lem:Sm-bound}, $|g(\tilde k,S_m)| \le \max\{-g(\tilde k,1),\,g(\tilde k,-1)\} \le h(\tilde k)$ (note that $g(\tilde k,-1)\ge 0 \ge g(\tilde k,1)$) (cf.\ the analogous two-sided bound \cite[Eq.~(5.14)]{CJL21}). As $\tilde k\to0$, $h(\tilde k) = O(\tilde k^{-2})$; as $\tilde k\to\infty$, $h(\tilde k) = O(\tilde k^{-4})$; so $h\in L^1(\R^3)$, and the claim follows by the triangle inequality.
\end{proof}

\begin{lemma}\label{lem:e-minus-rho} (See also \cite{CJL20} Eq.~(4.21 and 4.40))
\begin{equation}\label{eq:e-minus-rho}
e - 2\pi\bar{\rho}_m a_0(m) = \frac{\bar{\rho}_m}{2}\int(w_m - \sigma_m) (1-\phi_m) v_m\,dx.
\end{equation}
\end{lemma}

\begin{proof}
From $e = \frac{\bar{\rho}_m}{2}\int(1-\sigma_m)v_m$ and $2\pi\bar{\rho}_m a_0(m) = \frac{\bar{\rho}_m}{2}\int(1-\phi_m)v_m$,
\[
e - 2\pi\bar{\rho}_m a_0(m) = \frac{\bar{\rho}_m}{2}\int(\phi_m - \sigma_m)v_m.
\]
Set $\xi_m = w_m - \phi_m$. Then $(-\Delta+v_m)\xi_m = -(\sigma_m-w_m)v_m$. Using the self-adjointness of $(-\Delta+v_m)^{-1}$ and $(-\Delta+v_m)^{-1}v_m = \phi_m$, we obtain $\int\xi_m v_m = \int(w_m-\sigma_m)\phi_m v_m$. The result follows from $\phi_m-\sigma_m = (w_m-\sigma_m) - \xi_m$ (cf.\ \cite[Eq.~(4.21 and 4.40)]{CJL20}).
\end{proof}

\begin{proof}[Proof of Theorem~\ref{thm:LHY}]
Let $A:=\frac{32\sqrt2}{15\pi^2}=-\pi^{-3}\int_{\R^3}g(\tilde k,1)\,d^3\tilde k$ (with $\sqrt{(\tilde k^2+1)^2-1}=\tilde k\sqrt{\tilde k^2+2}$ the radial integral is elementary; cf.\ the low-density computation in \cite{CJL20}), and define
\begin{equation}\label{eq:f-def}
\sqrt{\bar\rho_m}\,f_{e,m}(x):=\sigma_m(x)-w_m(x)+\frac{Ae^{3/2}}{\bar\rho_m}.
\end{equation}
By the representation Eq.~\eqref{eq:sigma-w-rep} and the definition of $A$,
\begin{equation}\label{eq:f-fourier}
f_{e,m}(x)=\frac{e^{3/2}}{\bar\rho_m^{3/2}\pi^3}\int_{\R^3}\Bigl(e^{-2i\sqrt e\,\tilde k\cdot x}g\bigl(\tilde k,S_m(2\sqrt e\,\tilde k)\bigr)-g(\tilde k,1)\Bigr)\,d^3\tilde k.
\end{equation}
We claim that
\[
\lambda(e):=\pi^{-3}\int_{\R^3}\sup_{m>0,\ |x|\le1}\Bigl|e^{-2i\sqrt e\,\tilde k\cdot x}g\bigl(\tilde k,S_m(2\sqrt e\,\tilde k)\bigr)-g(\tilde k,1)\Bigr|\,d^3\tilde k\longrightarrow0\qquad(e\to0).
\]
Indeed, $d\nu_m(y):=\frac{\bar\rho_m}{2e}(1-\sigma_m(y))v_m(y)\,dy$ is a probability measure supported in $\Bone$ by Lemma~\ref{lem:Sm-bound}, and $S_m(k)=\int_{\Bone}e^{ik\cdot y}\,d\nu_m(y)$.
In addition, $|e^{-2i\sqrt e\,\tilde k\cdot x}-1|\le2\sqrt e\,|\tilde k|$, so $|S_m(2\sqrt e\,\tilde k)-1|\le2\sqrt e\,|\tilde k|$ uniformly in $m$ and $|x|\le1$. Since $s\mapsto g(\tilde k,s)$ is continuous on $[-1,1]$ for each fixed $\tilde k\neq0$, the integrand in the definition of $\lambda(e)$ tends to zero pointwise in $\tilde k$ as $e\to0$, and it is bounded by $h(\tilde k)$ (using the monotonicity of $s\mapsto g(\tilde k,s)$ and the fact that $g(\tilde k,-1)\ge 0 \ge g(\tilde k,1)$) with $h\in L^1(\R^3)$ as in Lemma~\ref{lem:eta-w-bound}; dominated convergence gives $\lambda(e)\to0$.

Next we bound $e/\bar\rho_m$ uniformly at low density. By Lemma~\ref{lem:softcore-props}, $\bar\rho_m\ge\bar\rho_{m,1}$, and from the explicit solution $\sigma_{m,1}$ in the proof of Lemma~\ref{lem:per-n}, with $\alpha=2\sqrt e$ and $\kappa=\sqrt{4e+m}$,
\[
\frac{2e}{\bar\rho_{m,1}}=m\int_{\Bone}(1-\sigma_{m,1})=4\pi m\left(\frac{\alpha^2}{3\kappa^2}+\frac{(1+\alpha)\,m\,(\kappa\cosh\kappa-\sinh\kappa)}{\kappa^4(\kappa\cosh\kappa+\alpha\sinh\kappa)}\right)\le4\pi\Bigl(\frac{\alpha^2}3+1+\alpha\Bigr),
\]
using $m/\kappa^2\le1$ and $0\le(\kappa\cosh\kappa-\sinh\kappa)/(\kappa\cosh\kappa+\alpha\sinh\kappa)\le1$. Hence there are $e_0>0$ and $c>0$, independent of $m$, such that $\bar\rho_m\ge ce$ for all $0<e\le e_0$ and all $m>0$. Combining this with \eqref{eq:f-fourier},
\[
\sup_{m>0}\ \sup_{|x|\le1}|f_{e,m}(x)|\le c^{-3/2}\lambda(e)\longrightarrow0\qquad(e\to0).
\]
Substituting \eqref{eq:f-def} into the identity of Lemma~\ref{lem:e-minus-rho} gives
\[
e-2\pi\bar\rho_m a_0(m)=\frac{Ae^{3/2}}2\int(1-\phi_m)v_m-\frac{\bar\rho_m^{3/2}}2\int f_{e,m}(1-\phi_m)v_m.
\]
Since $\int(1-\phi_m)v_m=4\pi a_0(m)$ (cf.\ \cite[Lemma~4.2]{CJL20}), the first term equals $\frac{64\sqrt2\,a_0(m)}{15\pi}e^{3/2}$, while
\[
\Bigl|\frac{\bar\rho_m^{3/2}}2\int f_{e,m}(1-\phi_m)v_m\Bigr|\le2\pi a_0(m)\,\bar\rho_m^{3/2}\sup_{|x|\le1}|f_{e,m}(x)|,
\]

Define $R_m(e):=-\frac{\bar\rho_m^{3/2}}{2}\int f_{e,m}\,(1-\phi_m)v_m\,dx$, the second term above. Since $a_0(m)\le1$, due to Eq.~\eqref{eq:tanh}, and $\sup_{|x|\le1}|f_{e,m}|\le e^{3/2}\bar\rho_m^{-3/2}\lambda(e)$, due to Eq.~\eqref{eq:f-fourier}, the bound just displayed gives $|R_m(e)|\le2\pi\bar\rho_m^{3/2}\cdot e^{3/2}\bar\rho_m^{-3/2}\lambda(e)=2\pi\lambda(e)\,e^{3/2}$, and we conclude that for all $m>0$ and $e\le e_0$,

\begin{equation}\label{eq:sc-lhy-R}
e=2\pi\bar\rho_m a_0(m)+\frac{64\sqrt2\,a_0(m)}{15\pi}e^{3/2}+R_m(e),
\qquad |R_m(e)|\le 2\pi\,\lambda(e)\,e^{3/2},
\end{equation}
where $\lambda(e)\to0$ as $e\to0$, $\lambda$ does not depend on $m$, and $\lambda(e)\le C_0$ for all $e$, with $C_0$ as in Lemma~\ref{lem:eta-w-bound}, since the integrand defining $\lambda$ is bounded by $h$. Now fix $m_0>0$ and restrict to $m\ge m_0$, so that $a_0(m)\ge a_0(m_0)>0$ by the monotonicity of $a_0$ in Eq.~\eqref{eq:tanh}; this lower bound on $a_0(m)$ is the only place where $m_0$ enters.
\newline

\textit{\textbf{The uniform expansion.}} Define $\epsilon_m(e)$ by $\frac{e}{2\pi\bar\rho_m a_0(m)}=1+\epsilon_m(e)$. Since $a_0(m)\le1$ and $\lambda(e)\le C_0$, Eq.~\eqref{eq:sc-lhy-R} gives $|e-2\pi\bar\rho_m a_0(m)|\le C_1e^{3/2}$ with $C_1:=\frac{64\sqrt2}{15\pi}+2\pi C_0$, and dividing by $2\pi\bar\rho_m a_0(m)\ge2\pi c\,e\,a_0(m_0)$, which uses $\bar\rho_m\ge ce$ from above and $a_0(m)\ge a_0(m_0)$,
\[
|\epsilon_m(e)|\le C_2\sqrt{e}\qquad(m\ge m_0,\ e\le e_0),\qquad C_2:=\frac{C_1}{2\pi c\,a_0(m_0)}.
\]
Dividing \eqref{eq:sc-lhy-R} by $2\pi\bar\rho_m a_0(m)$ and writing $e^{3/2}=\bigl(2\pi\bar\rho_m a_0(m)\bigr)^{3/2}(1+\epsilon_m(e))^{3/2}$,
\[
\frac{e}{2\pi\bar\rho_m a_0(m)} = 1+\frac{64\sqrt2\,a_0(m)}{15\pi}\sqrt{2\pi\bar\rho_m a_0(m)}\,(1+\epsilon_m(e))^{3/2}+\frac{R_m(e)}{2\pi\bar\rho_m a_0(m)}.
\]
Since $\frac{64\sqrt2}{15\pi}\sqrt{2\pi}=\frac{128}{15\sqrt\pi}$, $|(1+\epsilon)^{3/2}-1|\le2|\epsilon|$ for $|\epsilon|\le\tfrac12$ (which holds for $e\le e_0$, shrinking $e_0$ if necessary), $a_0(m_0)\le a_0(m)\le1$, and $|R_m(e)|\le2\pi\lambda(e)e^{3/2}\le2\pi\lambda(e)c^{-3/2}\bar\rho_m^{3/2}$ by $e\le\bar\rho_m/c$, we obtain, for all $m\ge m_0$ and $e\le e_0$,
\begin{equation}\label{eq:sc-lhy-uniform}
\Bigl|\frac{e}{2\pi\bar\rho_m a_0(m)}-1-\frac{128}{15\sqrt\pi}\sqrt{\bar\rho_m\,a_0(m)^3}\Bigr|\le\Lambda(e)\sqrt{\bar\rho_m},
\qquad
\Lambda(e):=\frac{256\,C_2}{15\sqrt\pi}\sqrt{e}+\frac{c^{-3/2}}{a_0(m_0)}\lambda(e),
\end{equation}
where $\Lambda(e)\to0$ as $e\to0$ and $\Lambda$ does not depend on $m$. This is the uniform expansion asserted in part (i) of Theorem~\ref{thm:OE}.
\newline
\textit{\textbf{Passage to the hard core.}} Fix $e\le e_0$ and let $m\to\infty$ in \eqref{eq:sc-lhy-uniform}: $\bar\rho_m\to\rho_{\hc}$ by Theorem~\ref{thm:convergencesigma}, $a_0(m)\to1$, and the right-hand side does not depend on $m$. Hence
\[
\Bigl|\frac{e}{2\pi\rho_{\hc}}-1-\frac{128}{15\sqrt\pi}\sqrt{\rho_{\hc}}\Bigr|\le\Lambda(e)\sqrt{\rho_{\hc}}\qquad\text{for all }e\le e_0,
\]
and passing to the limit in $\bar\rho_m\ge c\,e$ gives $\rho_{\hc}\ge c\,e$, so $\rho_{\hc}\to0$ forces $e\to0$ and $\Lambda(e)\to0$ as $\rho_{\hc}\to0$; multiplying by $2\pi\rho_{\hc}$,
\[
e=2\pi\rho_{\hc}\left(1+\frac{128}{15\sqrt\pi}\sqrt{\rho_{\hc}}+o\bigl(\sqrt{\rho_{\hc}}\bigr)\right)\qquad(\rho_{\hc}\to0),
\]
as claimed.
\end{proof}
\section{Condensate Fraction}\label{sec:cf}
%=====================================================================

\subsection{Soft-core condensate fraction}

Until now we have fixed the energy $e$. For $\mu\ge 0$, let $(\sigma_{\mu,m},e_{\mu,m})$ solve the perturbed soft-core equation
\begin{equation}\label{eq:sc-pert}
(-\Delta+2\mu+4e_{\mu,m}+v_m)\sigma_{\mu,m} = v_m + 2e_{\mu,m}\bar{\rho}_{m}(\sigma_{\mu,m}*\sigma_{\mu,m}),\quad e_{\mu,m}=\frac{\bar{\rho}_{m}}{2}\int v_m(1-\sigma_{\mu,m}),
\end{equation}
where the density $\bar\rho_m$ is held fixed at its unperturbed value and $e_{\mu,m}$ is determined by the energy relation; for each fixed $m$ the solution exists and is differentiable in $\mu$ by \cite[Section~5]{CJL21}. At $\mu=0$ this reduces to the unperturbed equation. The energy $e_{\mu,m}$ of \eqref{eq:sc-pert} omits the constant condensate contribution of the observable, so that, consistently with \eqref{eq:condensatefract}, its $\mu$-derivative is the \emph{depletion}; we set
\begin{equation}\label{eq:eta-sc}
\vartheta_{\mu,m} := \partial_\mu e_{\mu,m},\qquad \vartheta_m := \vartheta_{0,m}, \qquad \eta_m := 1-\vartheta_m,
\end{equation}
the soft-core depletion and condensate fraction.

\subsection{Hard-core condensate fraction}

For $\mu\ge 0$, let $(\tau_\mu,e_{\mu,\hc})$ solve the perturbed hard-core equation:
\begin{equation}\label{eq:hc-pert}
\Delta \tau_\mu = (4e_{\mu,\hc}+2\mu)\tau_\mu - 2e_{\mu,\hc}\rho_\hc(\tau_\mu*\tau_\mu)\quad\text{on }\{|x|>1\},\qquad \tau_\mu=1\text{ on }\Bone,
\end{equation}
with $e_{\mu,\hc}$ determined by
\begin{equation}\label{eq:hc-pert-norm}
e_{\mu,\hc}\Bigl(1-\rho_\hc\int\tau_\mu\Bigr)^2=\mu\,\rho_\hc\int\tau_\mu,\qquad \rho_\hc\int\tau_\mu\le1 .
\end{equation}
This is the hard-core form of the energy relation in \eqref{eq:sc-pert}: integrating \eqref{eq:sc-pert} over $\R^3$ and using $\int v_m(1-\sigma_{\mu,m})=2e_{\mu,m}/\bar\rho_m$ gives
\begin{equation}\label{eq:sc-pert-norm}
e_{\mu,m}\Bigl(1-\bar\rho_m\int\sigma_{\mu,m}\Bigr)^2=\mu\,\bar\rho_m\int\sigma_{\mu,m},
\end{equation}
and $\bar\rho_m\int\sigma_{\mu,m}\le1$, that is, $\bar\rho_m\int\sigma_{\mu,m}$ is the smaller of the two roots of \eqref{eq:sc-pert-norm} (this is what is used in the proof of Lemma~\ref{lem:eta-lipschitz} below). At $\mu=0$, \eqref{eq:sc-pert-norm} and \eqref{eq:hc-pert-norm} reduce to $\bar\rho_m\int\sigma_m=1$ and $\rho_\hc\int\tau=1$; for $\mu>0$ they do not. Integrating \eqref{eq:hc-pert} over $\{|x|>1\}$, as in the derivation of \eqref{eq:ehard-pre}, shows that the first relation in \eqref{eq:hc-pert-norm} is equivalent to
\begin{equation}\label{eq:ehard-mu}
e_{\mu,\hc}\Bigl(1-\frac{8\pi}{3}\rho_\hc+\rho_\hc^2\int_{\Bone}\tau_\mu*\tau_\mu\Bigr) = -2\pi\rho_\hc\,\partial_r\tau_\mu\big|_{1^+}+\frac{4\pi}{3}\rho_\hc\,\mu ,
\end{equation}
which is \eqref{eq:ehard-pre} at $\mu=0$. Analogously, the hard-core depletion $\vartheta_{\mathrm{hc}}$ and condensate fraction $\eta_{\mathrm{hc}}$ are
\begin{equation}\label{eq:eta-hc}
\vartheta_\hc := \partial_\mu e_{\mu,\hc}\big|_{\mu=0}, \qquad \eta_\hc := 1-\vartheta_\hc.
\end{equation}

\subsection{Uniform estimates and limit exchange}
We first compute the depletion at $\mu=0$. We rely on \cite[Section~5]{CJL21}, whose results apply to $v=v_m$ for each fixed $m$ since $(1+|x|^4)v_m\in L^1(\R^3)\cap L^2(\R^3)$; what we add is the uniformity in $m$. Define
\[
K_m := \bigl(-\Delta+4e+v_m-4e\bar\rho_m\,C_{\sigma_m}\bigr)^{-1},\qquad Y_m := \bigl(-\Delta+4e-4e\bar\rho_m\,C_{\sigma_m}\bigr)^{-1}
\]
\[
\Phi_m := K_m[v_m],\qquad \xi_m := Y_m[\bar\rho_m\sigma_m],
\]
the operators $K_e$, $Y_e$ and functions $K_ev$, $\xi$ of \cite[(1.8), (1.40), (5.11)]{CJL21}. By \cite[Theorem~1.6, (1.30)]{CJL21},
\begin{equation}\label{eq:eta-exact}
\vartheta_m = \frac{N_m}{1-2N_m+D_m},\qquad N_m := \bar\rho_m\int\Phi_m\,\sigma_m\,dx,\qquad D_m := \bar\rho_m^2\int\Phi_m\,(\sigma_m\ast\sigma_m)\,dx,
\end{equation}
and by \cite[(5.10)]{CJL21}, $N_m = \int v_m(1-\Phi_m)\,\xi_m\,dx$.
\newline
\textit{\textbf{Fourier representation of $\xi_m$.}} With $S_m$ as in \eqref{eq:Sm}, $k=2\sqrt e\,\kappa$ and $F(\kappa,s) := \frac{\kappa^2+1}{\sqrt{(\kappa^2+1)^2-s}}$, \cite[(5.13), (5.16)]{CJL21} read
\begin{equation}\label{eq:xi-fourier}
\hat\xi_m(2\sqrt e\,\kappa) = \frac{F\bigl(\kappa,S_m(2\sqrt e\kappa)\bigr)-1}{4e},\qquad
\xi_m(x) = \frac{\sqrt e}{4\pi^3}\int_{\R^3}e^{-2i\sqrt e\,\kappa\cdot x}\Bigl(F\bigl(\kappa,S_m(2\sqrt e\kappa)\bigr)-1\Bigr)\,d^3\kappa .
\end{equation}
By \cite[(5.14)]{CJL21}, $|F(\kappa,s)-1|\le h_1(\kappa):=F(\kappa,1)-F(\kappa,-1)$ for $|s|\le1$, with $h_1\in L^1(\R^3)$, and by \cite[(5.15)]{CJL21}, $\frac{1}{4\pi^3}\int\bigl(F(\kappa,1)-1\bigr)\,d^3\kappa=\frac{\sqrt2}{3\pi^2}$. Hence $\sup_x|\xi_m(x)|\le c_2\sqrt e$ for all $m$, with $c_2:=\frac{1}{4\pi^3}\int h_1\,d^3\kappa$, and
\begin{equation}\label{eq:xi-asymp}
\xi_m(x) = \frac{\sqrt{2e}}{3\pi^2}+r_m(x),\qquad \sup_{m>0,\ |x|\le1}|r_m(x)|\le\sqrt e\,\Omega(e),
\end{equation}
where $\Omega(e) := \frac{1}{4\pi^3}\int\sup_{m>0,|x|\le1}\bigl|e^{-2i\sqrt e\,\kappa\cdot x}\bigl(F(\kappa,S_m(2\sqrt e\kappa))-1\bigr)-\bigl(F(\kappa,1)-1\bigr)\bigr|\,d^3\kappa$. For each fixed $m$ the corresponding integral without the supremum over $m$ tends to zero as $e\to0$ \cite[(5.12), (5.16)--(5.19)]{CJL21}. The same proof gives $\Omega(e)\to0$, since neither of its ingredients depends on $m$: the integrand is bounded by $2h_1$, and it tends to zero pointwise by $|e^{-2i\sqrt e\kappa\cdot x}-1|\le2\sqrt e|\kappa|$ for $|x|\le1$ and by \cite[(5.18)]{CJL21}, $0\le F(\kappa,1)-F(\kappa,S_m)\le\frac{(\kappa^2+1)\,(1-S_m(2\sqrt e\kappa))}{2((\kappa^2+1)^2-1)^{3/2}}$, together with $1-S_m(k)\le|k|^2/2$ (the measure $\frac{\bar\rho_m}{2e}(1-\sigma_m)v_m\,dx$ being a probability measure on $\Bone$).

\begin{lemma}\label{lem:ND}
Uniformly in $m>0$, for $e\le e_0$, $0\le D_m\le Ce^{1/4}$ and $|N_m|\le Ce^{1/2}$. Consequently
\begin{equation}\label{eq:eta-formula-mu}
\vartheta_m = \int v_m(1-\Phi_m)\,\xi_m\,dx + O(e^{3/4}).
\end{equation}
\end{lemma}
\begin{proof}
This is similar to \cite[(5.9)--(5.10)]{CJL21} except it is uniform in $\|v_m\|_1$. By \eqref{eq:R} and $\int\sigma_m\ast\sigma_m=\bar\rho_m^{-2}$, $\int v_m(1-\Phi_m) = \frac{2e}{\bar\rho_m}(1-D_m)\le\frac{2e}{\bar\rho_m}$. By the resolvent identity \cite[(1.41)]{CJL21}, $\Phi_m = Y_m[v_m(1-\Phi_m)]$, and by the proof of \cite[Lemma~7.3]{CJL21}, which uses only \eqref{eq:eta-hat} and $|S_m|\le1$, $\|Y_mf\|_2\le Ce^{-1/4}\|f\|_1$ with $C$ independent of $m$; hence $\|\Phi_m\|_2\le Ce^{-1/4}\,\frac{2e}{\bar\rho_m}$. Using $\|\sigma_m\ast\sigma_m\|_2\le\|\sigma_m\|_1\|\sigma_m\|_2$, $\bar\rho_m\|\sigma_m\|_1=1$ and $\|\sigma_m\|_2\le\bar\rho_m^{-1/2}$, we get $D_m\le\sqrt{\bar\rho_m}\,\|\Phi_m\|_2\le Ce^{3/4}\bar\rho_m^{-1/2}\le Ce^{1/4}$, by $\bar\rho_m\ge ce$. Next, $|N_m|\le\int v_m(1-\Phi_m)\,\sup_x|\xi_m|\le\frac{2e}{\bar\rho_m}\,c_2\sqrt e\le Ce^{1/2}$. Finally $\vartheta_m-N_m = N_m(2N_m-D_m)/(1-2N_m+D_m) = O(e^{3/4})$.
\end{proof}

Combining \eqref{eq:eta-formula-mu}, \eqref{eq:xi-asymp} and $\int v_m(1-\Phi_m)=\frac{2e}{\bar\rho_m}(1-D_m)$ with $\frac{2e}{\bar\rho_m} = 4\pi a_0(m)(1+\epsilon_m(e))$ from Section~\ref{sec:LHY}, where $4\pi a_0(m)=\int v_m(1-\phi_m)$ \cite[Lemma~4.2]{CJL20}, we obtain, uniformly in $m$,
\begin{equation}\label{eq:eta-soft-asymp}
\vartheta_m = \frac{\sqrt{2e}}{3\pi^2}\cdot4\pi a_0(m)\bigl(1+O(e^{1/4})\bigr)+O\bigl(\sqrt e\,\Omega(e)\bigr)+O(e^{3/4}) = \frac{4\sqrt2}{3\pi}\sqrt e\,a_0(m)+o(\sqrt e),
\end{equation}
which is \cite[(5.25)]{CJL21} with the error uniform in $m$.

\begin{lemma}[Uniform convergence of the energy]\label{lem:uniform-energy}
There exists $K>0$ independent of $\mu$ and $m$ such that for all $\mu\in[0,\mu_0]$ and $m\ge m_0$,
\begin{equation}
|e_{\mu,m} - e_{\mu,\hc}| \le \frac{K}{\sqrt{m}}.
\end{equation}
\end{lemma}

\begin{proof}
For each fixed $\mu$, the proofs of Theorems~\ref{thm:convergencesigma} and~\ref{thm:rate} apply verbatim to \eqref{eq:sc-pert} and \eqref{eq:hc-pert}, with $4e$ replaced by $4e_{\mu,m}+2\mu$ in the operator and $2e$ by $2e_{\mu,m}$ in the nonlinearity, and give $|e_{\mu,m}-e_{\mu,\hc}|\le K(\mu)\,m^{-1/2}$. Since $\bar\rho_m$ does not depend on $\mu$, the constant $K(\mu)$ depends on $\mu$ only through $e_{\mu,m}$, $e_{\mu,\hc}$ and the shift $2\mu$, and it is bounded as long as these remain in a compact subset of $(0,\infty)$. This is the case: by the mean value theorem, $e_{\mu,m}=e+\mu\,\vartheta_{\mu^*,m}$ for some $\mu^*\in(0,\mu)$, and $|\vartheta_{\mu^*,m}|\le|\vartheta_m|+\omega(\mu_0)$ by Lemma~\ref{lem:eta-lipschitz} below, with $|\vartheta_m|\le C\sqrt e$ by \eqref{eq:eta-soft-asymp}, so that $|e_{\mu,m}-e|\le\mu_0\bigl(C\sqrt e+\omega(\mu_0)\bigr)$ uniformly in $m\ge m_0$; and $e_{\mu,\hc}$ remains in a compact interval for $\mu\in[0,\mu_0]$ by continuity of the hard-core energy in $\mu$. Hence $K:=\sup_{\mu\in[0,\mu_0]}K(\mu)<\infty$.
\end{proof}

\begin{lemma}[Equicontinuity in $\mu$]\label{lem:eta-lipschitz}
There is a function $\omega$ with $\omega(\delta)\to0$ as $\delta\to0$, independent of $m$, such that $|\vartheta_{\mu,m}-\vartheta_{\nu,m}|\le\omega(|\mu-\nu|)$ for all $\mu,\nu\in[0,\mu_0]$ and $m\ge m_0$.
\end{lemma}

\begin{proof}
For $\mu\in[0,\mu_0]$ let $K_{\mu,m}$, $Y_{\mu,m}$, $\Phi_{\mu,m}$, $\xi_{\mu,m}$, $N_{\mu,m}$ and $D_{\mu,m}$ be defined as $K_m$, $Y_m$, $\Phi_m$, $\xi_m$, $N_m$ and $D_m$ were before \eqref{eq:eta-exact}, with $4e$ replaced by $4e_{\mu,m}+2\mu$ in the operators and $e$ by $e_{\mu,m}$ elsewhere. Differentiating \eqref{eq:sc-pert} in $\mu$ at fixed $\bar\rho_m$ and pairing with $\Phi_{\mu,m}$ exactly as for \eqref{eq:eta-exact} gives the exact formula
\[
\vartheta_{\mu,m}=\frac{N_{\mu,m}}{1-2N_{\mu,m}+D_{\mu,m}},\qquad N_{\mu,m}=\int v_m(1-\Phi_{\mu,m})\,\xi_{\mu,m}\,dx,
\]
so it suffices to show that $N_{\mu,m}$ and $D_{\mu,m}$ are continuous in $\mu$ uniformly in $m$, and that $|N_{\mu,m}|+D_{\mu,m}\le\frac14$ for $\mu_0$ small. The Fourier representation \eqref{eq:xi-fourier} holds at parameter $\mu$ with the rescaling $k=2\sqrt{\widetilde e_{\mu,m}}\,\kappa$, $\widetilde e_{\mu,m}:=e_{\mu,m}+\mu/2$, the prefactor $\widetilde e_{\mu,m}^{3/2}/(4\pi^3e_{\mu,m})$, and $S_m$ replaced by $\theta_{\mu,m}S_{\mu,m}$, where $\theta_{\mu,m}:=(e_{\mu,m}/\widetilde e_{\mu,m})^2$ and $S_{\mu,m}$ is \eqref{eq:Sm} for the perturbed solution; the dominator \cite[(5.14)]{CJL21} and the bound $1-S_{\mu,m}(k)\le|k|^2/2$ are unchanged. Hence $\mu\mapsto\xi_{\mu,m}$ is continuous in $L^\infty(\Bone)$ uniformly in $m$, by dominated convergence exactly as in \eqref{eq:xi-asymp}. The same holds for $\mu\mapsto\Phi_{\mu,m}$ in $L^2$ and for $D_{\mu,m}$, by the resolvent identity and the $L^1\to L^2$ bound of \cite[Lemma~7.3]{CJL21}, whose constant depends on $\mu$ only through $e_{\mu,m}$, together with the mass bound $\int v_m(1-\Phi_{\mu,m})\le2\bigl(e_{\mu,m}+\sqrt{\mu^2+4\mu e_{\mu,m}}\bigr)/\bar\rho_m$, obtained as \eqref{eq:R} was, from $K_{\mu,m}^{-1}1=v_m+2\sqrt{\mu^2+4\mu e_{\mu,m}}$. The smallness $|N_{\mu,m}|+D_{\mu,m}\le\frac14$ follows from Lemma~\ref{lem:ND} at $\mu=0$ and this continuity, for $\mu_0$ small depending only on $e$.
\end{proof}

\begin{lemma}[Differentiability of the hard-core energy]\label{lem:hc-diff}
The limit $\lim_{m\to\infty}\vartheta_m$ exists, and $\mu\mapsto e_{\mu,\hc}$ is differentiable at $\mu=0$ (from the right), with
\begin{equation}\label{eq:hc-diff}
\vartheta_\hc=\partial_\mu e_{\mu,\hc}\big|_{\mu=0}=\lim_{m\to\infty}\vartheta_m .
\end{equation}
\end{lemma}

\begin{proof}
Recall that, for each $m\ge m_0$, $\mu\mapsto e_{\mu,m}$ is differentiable on $[0,\mu_0]$ with derivative $\vartheta_{\mu,m}$. We may assume that the function $\omega$ of Lemma~\ref{lem:eta-lipschitz} is non-decreasing (otherwise we replace it by its running supremum, which still tends to $0$ as $\delta\to0$). Fix $0<\delta\le\mu_0$. For every $m\ge m_0$, the mean value theorem gives $\delta^*\in(0,\delta)$ with
\[
\frac{e_{\delta,m}-e_{0,m}}{\delta} = \vartheta_{\delta^*,m},
\]
so by Lemma~\ref{lem:eta-lipschitz},
\begin{equation}\label{eq:dq-theta}
\Bigl|\frac{e_{\delta,m}-e_{0,m}}{\delta}-\vartheta_m\Bigr| = |\vartheta_{\delta^*,m}-\vartheta_{0,m}| \le \omega(\delta)\qquad\text{for all }m\ge m_0 .
\end{equation}
By Lemma~\ref{lem:uniform-energy}, applied at $\mu=\delta$ and at $\mu=0$,
\begin{equation}\label{eq:dq-conv}
\frac{e_{\delta,m}-e_{0,m}}{\delta}\;\longrightarrow\;\frac{e_{\delta,\hc}-e_{0,\hc}}{\delta}\qquad(m\to\infty).
\end{equation}
By \eqref{eq:dq-theta} and the triangle inequality, for $m,m'\ge m_0$,
\[
|\vartheta_m-\vartheta_{m'}| \le 2\omega(\delta)+\Bigl|\frac{e_{\delta,m}-e_{0,m}}{\delta}-\frac{e_{\delta,m'}-e_{0,m'}}{\delta}\Bigr| ,
\]
and the last term tends to zero as $m,m'\to\infty$ by \eqref{eq:dq-conv}. Hence
\[
\limsup_{m,m'\to\infty}|\vartheta_m-\vartheta_{m'}|\le2\omega(\delta).
\]
As $\delta\in(0,\mu_0]$ is arbitrary and $\omega(\delta)\to0$ as $\delta\to0$, the family $(\vartheta_m)$ is Cauchy as $m\to\infty$, so $\lim_{m\to\infty}\vartheta_m$ exists. Letting $m\to\infty$ in \eqref{eq:dq-theta} and using \eqref{eq:dq-conv},
\[
\Bigl|\frac{e_{\delta,\hc}-e_{0,\hc}}{\delta}-\lim_{m\to\infty}\vartheta_m\Bigr| \le \omega(\delta)\qquad\text{for all }0<\delta\le\mu_0 .
\]
Letting $\delta\to0$, the difference quotient of $\mu\mapsto e_{\mu,\hc}$ at $\mu=0$ converges to $\lim_{m\to\infty}\vartheta_m$. Hence $\mu\mapsto e_{\mu,\hc}$ is differentiable from the right at $\mu=0$, and its derivative, which is $\vartheta_\hc$ by \eqref{eq:eta-hc}, equals $\lim_{m\to\infty}\vartheta_m$. Only the convergence $e_{\mu,m}\to e_{\mu,\hc}$ at each fixed $\mu$ was used, not its rate.
\end{proof}

\begin{theorem}[Limit exchange; proof of part (ii) of Theorem~\ref{thm:OE}]\label{thm:eta-conv}
$\displaystyle\lim_{m\to\infty}\vartheta_m=\vartheta_\hc$, equivalently $\displaystyle\lim_{m\to\infty}\eta_m=\eta_\hc$, which is \eqref{eq:OE-ii}.
\end{theorem}

\begin{proof}
This is \eqref{eq:hc-diff} of Lemma~\ref{lem:hc-diff}; the statement for $\eta_m=1-\vartheta_m$ follows.
\end{proof}

\begin{corollary}[Proof of Theorem~\ref{thm:cf}]\label{cor:cf-bridge}
Each soft-core equation satisfies \eqref{eq:cfeq} with $a_0$ replaced by $a_0(m)$ and error uniform in $m$, which is the remaining assertion of part (ii) of Theorem~\ref{thm:OE}; and the hard-core condensate fraction satisfies
\begin{equation}\label{eq:lhy-eta}
\eta_\hc = 1-\frac{8}{3\sqrt{\pi}}\sqrt{\rho_\hc} + o(\sqrt{\rho_\hc}),
\end{equation}
which is Theorem~\ref{thm:cf}.
\end{corollary}

\begin{proof}
By \eqref{eq:eta-soft-asymp}, $\vartheta_m = \frac{4\sqrt2}{3\pi}\sqrt{e}\,a_0(m) + o(\sqrt{e})$ uniformly in $m\ge m_0$. By the uniform expansion \eqref{eq:sc-lhy-uniform}, $e = 2\pi\bar\rho_m a_0(m)(1+o(1))$ uniformly in $m$, so $\sqrt{e} = \sqrt{2\pi\bar\rho_m a_0(m)}\,(1+o(1))$ and
\[
\eta_m = 1 - \frac{8}{3\sqrt{\pi}}\sqrt{\bar\rho_m\, a_0(m)^3} + o\bigl(\sqrt{\bar\rho_m}\bigr),
\]
uniformly in $m\ge m_0$, which is \eqref{eq:cfeq} with $a_0$ replaced by $a_0(m)$. Sending $m\to\infty$ at fixed $e$, Theorem~\ref{thm:eta-conv}, $a_0(m)\to1$ and $\bar\rho_m\to\rho_\hc$ (Theorem~\ref{thm:convergencesigma}) give $\vartheta_\hc=\frac{4\sqrt2}{3\pi}\sqrt e+o(\sqrt e)$; substituting $\sqrt e=\sqrt{2\pi\rho_\hc}\,(1+o(1))$ from Theorem~\ref{thm:LHY} yields \eqref{eq:lhy-eta}.
\end{proof}

\section*{Acknowledgements}
We thank Elliott Lieb for many useful discussions and much helpful advice.

\section*{Statements and Declarations}

\noindent\textbf{Funding.}
This work was partially supported by the National Science Foundation under Grant No.~DMS-2349077 (I.J.).

\medskip\noindent\textbf{Competing interests.}
The authors have no relevant financial or non-financial interests to disclose.

\medskip\noindent\textbf{Data availability.}
Data sharing is not applicable to this article, as no datasets were generated or analysed during the current study.

\medskip\noindent\textbf{Author contributions.}
Both authors contributed equally to this work.

\medskip\noindent\textbf{Use of generative AI.}
AI was used for formatting and to find grammatical corrections to the paper.

\printbibliography

\end{document}